\documentclass{amsart}
 \usepackage[foot]{amsaddr}
\usepackage{enumitem}
\usepackage{amsmath}
\usepackage{amssymb}
\usepackage{centernot}
\usepackage{mathtools}
\usepackage{xcolor}

\makeatletter
\@namedef{subjclassname@2020}{\textup{2020} Mathematics Subject Classification}
\makeatother

\newcommand{\pres}[3]{\textnormal{#1} \langle #2 \mid #3 \rangle}

\newcommand{\lr}[1]{\xleftrightarrow{}_{#1}}
\newcommand{\lra}[1]{\xleftrightarrow{\ast}_{#1}}
\newcommand{\xra}[1]{\xrightarrow{\ast}_{#1}}
\newcommand{\xr}[1]{\xrightarrow{}_{#1}}

\newcommand{\hs}{\heartsuit}
\newcommand{\ths}{\widetilde{\heartsuit}}
\newcommand{\trev}{\text{rev}}
\newcommand{\dD}{\overline{\Delta}}
\newcommand{\CF}{\mathcal{C}_{\operatorname{cf}}}

\newcommand{\DCF}{\mathcal{C}_{\operatorname{dcf}}}
\newcommand{\REG}{\mathcal{C}_{\operatorname{reg}}}
\newcommand{\IND}{\mathcal{C}_{\operatorname{ind}}}

\newcommand{\cc}{\mathcal{C}}

\newcommand{\cR}{\mathcal{R}}
\newcommand{\cl}{\mathcal{L}}

\DeclareMathOperator{\AFL}{AFL}
\DeclareMathOperator{\Rep}{Rep}
\newcommand{\DRep}{\Rep\Delta}

\DeclareMathOperator{\IP}{IP}
\DeclareMathOperator{\WP}{WP}
\DeclareMathOperator{\InvP}{InvP}

\newtheorem{theorem}{Theorem}[section] 
\newtheorem*{theorems}{Theorem}
\newtheorem{lemma}[theorem]{Lemma}     
\newtheorem{corollary}[theorem]{Corollary}
\newtheorem*{corollarys}{Corollary}
\newtheorem{proposition}[theorem]{Proposition}
\newtheorem{theoremB}{Theorem}

\theoremstyle{definition}
\newtheorem{question}{Question}
\newtheorem{example}{Example}
\newtheorem{remark}{Remark}
\newtheorem{conjecture}{Conjecture}

\numberwithin{equation}{section}

\title[On the Word Problem for Special Monoids]
 {On the Word Problem for Special Monoids} 



\author{Carl-Fredrik Nyberg-Brodda}
\address{Department of Mathematics, Alan Turing Building, University of Manchester, UK.}
\email{carl-fredrik.nybergbrodda@manchester.ac.uk}
\thanks{The author is currently a research associate, funded by the Dame Kathleen Ollerenshaw Trust, at the University of Manchester, United Kingdom}

\date{today}

\date{\today}

\begin{document}

\begin{abstract}
A monoid is called \textit{special} if it admits a presentation in which all defining relations are of the form $w = 1$. Every group is special, but not every monoid is special. In this article, we describe the language-theoretic properties of the word problem, in the sense of Duncan \& Gilman, for special monoids in terms of their group of units. We prove that a special monoid has context-free word problem if and only if its group of units is virtually free, giving a full generalisation of the Muller-Schupp theorem. This fully answers, for the class of special monoids, a question posed by Duncan \& Gilman in 2004. We describe the congruence classes of words in a special monoid, and prove that these have the same language-theoretic properties as the word problem. This answers a question first posed by Zhang in 1992. As a corollary, we prove that it is decidable (in polynomial time) whether a special one-relation monoid has context-free word problem. This completely answers another question from 1992, also posed by Zhang.
\end{abstract}

\maketitle

\noindent The study of algorithmic properties of algebraic structures has a rich and fascinating history, linking with a vast number of different areas of mathematics. One of the central problems for an algebraic structure is the \textit{word problem}. This is the problem of, given a presentation of the structure in terms of generators and relations, deciding whether or not two given words over the generators represent the same element of the structure. This problem can be dated back to the early 20th century in the investigations by Dehn and Thue \cite{Dehn1911, Thue1914}. Natural algebraic structures to investigate with respect to this problem are monoids, i.e. semigroups with an identity element. A \textit{special} monoid is one which admits a presentation in which all the defining relations are of the form $w = 1$. Special monoids were first defined, and given their name, in 1958 by Tse{\u{\i}}tin \cite[p. 178]{Tseitin1958}, although special monoids with a single defining relation $w=1$ had implicitly already been studied in 1914 by Thue \cite[\S3, \S5--8]{Thue1914}. Indeed, \cite[Problem~II]{Thue1914} asks precisely for a solution to the word problem for special one-relation monoids $\pres{Mon}{A}{w=1}$, and the remainder of Thue's article deals with solving special cases of precisely this problem.\footnote{Thue, ibid. p. 8: ``Man kann sich nun die gro{\ss}e Aufgabe (II) stellen: [...] Das eigentliche Ziel unserer Abhandlung besteht nur darin, die L\"osung einiger Beispiele dieser Aufgabe zu geben''.} Thus special monoids form one of the cornerstones of combinatorial semigroup theory. 

 The first structured approach to special monoids, however, was due to Adian and his student Makanin in the 1960s \cite{Adian1960, Adian1962, Adian1966, Makanin1966, Makanin1966b}. The usage of the term \textit{special} for such monoids was subsequently exported beyond the Soviet Union; it appears to have been first introduced to the French school by Lallement \cite[p. 371]{Lallement1974}, and to the German by Jantzen \cite{Jantzen1981}, both of whom explicitly borrow the phrase from Adian \cite{Adian1966}. Other terms in use, before \textit{special} became standard, were \textit{unitary}, used by e.g. Book \cite{Book1982}, and \textit{trivial}, used by e.g. Cochet \cite{Cochet1971}; the usage of this latter term is (justifiably) derided by Jantzen \cite[p. 72]{Jantzen1981}. A rewriting of Adian's and Makanin's proofs and results into the language of rewriting systems was later made by Zhang in the early 1990s. Other authors have also investigated special monoids, see e.g. \cite{Benois1973, Kashintsev1978, Kashintsev1993}.

Arguably the most important result to date on special monoids is the following, sometimes referred to as the Adian-Makanin Theorem: the group of units of a $k$-relation special monoid is a $k$-relator group, and a finitely presented special monoid has decidable word problem if and only if its group of units does \cite{Adian1960, Adian1966, Makanin1966}. An example of an immediate corollary to this result is that the word problem for any special monoid with a single defining relation $w=1$ is decidable, using Magnus' solution to the word problem for one-relator groups \cite{Magnus1932}. Hence, given a special monoid $M$, the group of units $U(M)$ plays a key r\^ole in understanding the monoid. We shall presently see this theme reinforced in this article.

From another direction, the methods of formal language theory have been highly successful in applications to group theory. This was initiated in 1971 by An{\={\i}}s{\={\i}}mov \cite{Anisimov1971}, and considers, for a group $G$ generated by a set $A$, the set of all words over $A$ which represent the identity element. The realisation that the language-theoretic properties of this set of words can reveal structural information about the group in question unlocked an entirely new angle of approach to group theory. The set is commonly referred to as the \textit{word problem} for the group, where the definite article is supported by the fact that the language-theoretic properties of this set do not generally depend on the generating set chosen for the group. An{\={\i}}s{\={\i}}mov showed that the word problem of a group is a regular language if and only if the group is finite.\footnote{However, An{\={\i}}s{\={\i}}mov remarks that this theorem can be immediately deduced already from results by Glu\v{s}kov on automatic partitions from 10 years prior \cite{Glushkov1961}.} Further, he showed that the class of groups with context-free word problem is closed under free products, and that it does not contain the group $\mathbb{Z}^2$. He later proved that every context-free group has the Howson property \cite{Anisimov1975}. A full and striking characterisation of the class of context-free groups was then provided: a finitely generated group has context-free word problem if and only if it is virtually free. This was proved in 1983 by Muller \& Schupp \cite{Muller1983} (when supplemented by a result by Dunwoody's \cite{Dunwoody1985}). This famous result -- the Muller-Schupp theorem -- is now the basis of a large literature and much active research, see e.g. \cite{Gilman1996, Gilman2007, Diekert2013, Araujo2017}.  

With the successes of language-theoretic techniques applied to the theory of groups evident, efforts were made to translate the relevant definitions and results to semigroup theory. In 2004, Duncan \& Gilman \cite{Duncan2004} defined a language which provides an elegant analogy for monoids of the word problem for groups, and which takes into account the fact that the set of words representing the identity element in a monoid need not be particularly enlightening. The language they defined, also referred to as the \textit{word problem}, enjoys the same types of invariance under choice of generating set as in the group case, and the analogy of An{\={\i}}s{\={\i}}mov's result remains true: a monoid has regular word problem if and only if it is finite \cite[p. 523]{Duncan2004}. In view of the Muller-Schupp theorem, a natural question is thus characterising which monoids have context-free word problem; this was the first question asked by Duncan \& Gilman about their definition of the word problem \cite[Question~4]{Duncan2004}. This remains a wide open problem.

The goal of this article is to investigate (and fully answer) Duncan \& Gilman's question in the setting of special monoids. More generally, we will ask:

\begin{question}\label{Quest: Alg struc of C-monoids?}
Let $\mathcal{C}$ be a class of languages. What is the algebraic structure of a special monoid with word problem in $\mathcal{C}$?
\end{question}

Duncan \& Gilman's question is thus the Question~\ref{Quest: Alg struc of C-monoids?} for $\cc = \CF$, the class of context-free languages. In this article, we will answer Question~\ref{Quest: Alg struc of C-monoids?} completely for $\cc = \CF$, but also for much more general classes $\cc$, namely when $\cc$ a super-$\AFL$ (in the sense of Greibach \cite{Greibach1970}) closed under reversal. Specifically, we have the following main theorem: 

\begin{theoremB}
Let $M$ be a finitely presented special monoid. Let $\cc$ be a super-$\AFL$ closed under reversal. Then $M$ has word problem in $\cc$ if and only if its group of units $U(M)$ has word problem in $\cc$.
\end{theoremB}

We deduce the following (Corollary~\ref{Cor:MS}) by applying the Muller-Schupp theorem:

\begin{corollarys}
Let $M$ be a finitely presented special monoid. Then $M$ has context-free word problem if and only if the group of units $U(M)$ is virtually free.
\end{corollarys}

As every group is a special monoid, this is a strict generalisation of the Muller-Schupp theorem for groups. Using Theorem~\ref{Thm:Main_thm}, we can deduce many language-theoretic properties of special monoids. First, we prove a result (Theorem~\ref{Thm:Equiv_of_three_notions}) which we do not state in full generality here, but which has as a special case:

\begin{theorems}
Let $M$ be a finitely presented special monoid. Then $M$ has context-free word problem if and only if the set of all words representing the identity element of $M$ is a context-free language. 
\end{theorems}

That is, the usual group-theoretic approach of considering the language of words equal to the identity element is also applicable to special monoids. Using this, we completely answer a question first posed in 1992 by Zhang. We also deduce some results regarding decision problems; for example, we prove (Theorem~\ref{Thm:RSMP_for_VF}) that for a special monoid with virtually free group of units, the rational subset membership problem is decidable. As another application, we find the second main theorem of this article:

\begin{theoremB}
Let $M = \pres{Mon}{A}{w=1}$ be a special one-relation monoid. Then it is decidable whether $M$ has context-free word problem. 
\end{theoremB}

In fact, this can be decided in polynomial time. Whether this was the case was asked by Zhang in 1992, and so Theorem~\ref{Thm:Dec_if_w=1_CF} answers this question. This theorem can be seen as a first step in a program to completely classify the language-theoretic properties of one-relation monoids with respect to various classes of languages. 

We also investigate Question~\ref{Quest: Alg struc of C-monoids?} for some classes of languages $\cc$ which are not super-$\AFL$s. Namely, we consider the case when $\cc$ is the class of regular languages or the class of deterministic context-free languages; we (easily) classify the special monoids with regular word problem, and prove that any special monoid with deterministic context-free word problem is right cancellative (Proposition~\ref{Prop:DCF=>RC+U(M)CF}). We end with a discussion of a property which a special monoid presentation may have (being \textit{infix}), and ask whether every special monoid admits an infix presentation. When the group of units is trivial, we answer this affirmatively (Proposition~\ref{Prop:UMtrivial=>Infix}).

\setcounter{theoremB}{0}

\clearpage
\section{Introduction}\label{Sec: Intro+pieces}

\noindent We assume the reader is familiar with formal language theory, including regular and context-free languages, as well as some elementary properties of codes. For some background on this, and other topics in formal language theory, we refer the reader to standard books on the subject \cite{Harrison1978, Hopcroft1979, Berstel1985}. The paper also assumes familiarity with the basics of the theory of monoid and group presentations, which will be written as $\pres{Mon}{A}{R}$ and $\pres{Gp}{A}{R}$, respectively. In particular, we assume the reader is familiar with the word problem as a decision problem. For further background information and examples of this theory, see e.g. \cite{Adian1966,Magnus1966,Neumann1967, Lyndon1977,Campbell1995}.

\subsection{Basic notation}

Let $A$ be a finite alphabet, and let $A^\ast$ denote the free monoid on $A$, with identity element denoted $\varepsilon$ or $1$, depending on the context. Let $A^+$ denote the free semigroup on $A$, i.e. $A^+ = A^\ast \setminus \{ \varepsilon\}$. For $u, v \in A^\ast$, by $u \equiv v$ we mean that $u$ and $v$ are the same word. For $w \in A^\ast$, we let $|w|$ denote the \textit{length} of $w$, i.e. the number of letters in $w$. We have $|\varepsilon| = 0$. If $w \equiv a_1 a_2 \cdots a_n$ for $a_i \in A$, then we let $w^\trev$ denote the \textit{reverse} of $w$, i.e. the word $a_n a_{n-1} \cdots a_1$. Note that $^\trev \colon A^\ast \to A^\ast$ is an anti-homomorphism, i.e. $(uv)^\trev \equiv v^\trev u^\trev$ for all $u, v \in A^\ast$. If $X \subseteq A^\ast$, then we let $X^\trev = \{ x^\trev \mid x \in X \}$. If the words $u, v \in A^\ast$ are equal in the monoid $M = \pres{Mon}{A}{R}$, then we denote this $u =_M v$. For a monoid $M$ generated by a finite set $A$, by which we mean there exists a surjective homomorphism $\pi \colon A^\ast \to M$, we define, for $w \in A^\ast$, the set
\[
\Rep_A^M(w) := \{ u \in A^\ast \mid u =_M w \}.
\]
This set $\Rep_A^M(w) \subseteq A^\ast$ is the set of \textit{representatives} of the element $\pi(w)$. Finally, for $X \subseteq A^\ast$ we let $\langle X \rangle$ denote the submonoid of $M$ generated by $X$, i.e. $\pi(X^\ast)$. 

We give some notation for rewriting systems. For an in-depth treatment and some terminology, see e.g. \cite{Book1982, Jantzen1988, Book1993}. A \textit{rewriting system} $R$ on $A$ is a subset of $A^\ast \times A^\ast$. An element of $R$ is called a \textit{rule}. If $(u, v) \in R$, we will sometimes denote this $(u \to v) \in R$. The system $R$ induces several relations on $A^\ast$. We will write $u \xr{R} v$ if there exist $x, y \in A^\ast$ and a rule $(\ell, r) \in T$ such that $u \equiv x\ell y$ and $v \equiv xry$. We let $\xra{R}$ denote the reflexive and transitive closure of $\xr{R}$. We denote by $\lra{R}$ the symmetric, reflexive, and transitive closure of $\xr{R}$. The relation $\lra{R}$ defines the least congruence on $A^\ast$ containing $R$. The monoid $\pres{Mon}{A}{R}$ is identified with the quotient $A^\ast / \lra{R}$. Let $u, v \in A^\ast$ and let $n \geq 0$. If there exist words $u_0, u_1, \dots, u_n \in A^\ast$ such that 
\[
u \equiv u_0 \xr{R} u_1 \xr{R} \cdots \xr{R} u_{n-1} \xr{R} u_n \equiv v,
\]
then we denote this $u \xr{R}^n v$, i.e. $u$ rewrites to $v$ in $n$ steps. Thus $\xra{R} = \bigcup_{n \geq 0} \xr{R}^n$.

A rewriting system $R \subseteq A^\ast \times A^\ast$ is said to be \textit{monadic} if $(u, v) \in R$ implies $|u| \geq |v|$ and $v \in A \cup \{ \varepsilon \}$. We say that $R$ is \textit{special} if $(u, v) \in R$ implies $v \equiv \varepsilon$. Every special system is monadic. A monadic rewriting system $R$ is said to be \textit{context-free} (resp. \textit{regular}) if for every $a \in A \cup \{ \varepsilon \}$, the language $\{ u \mid (u, a) \in R \}$ is context-free (resp. regular).

Let $G$ be a group with finite (group) generating set $A$, with $A^{-1}$ the set of inverses of the generators $A$. The language 
\[
\IP_{A \cup A^{-1}}^G := \{ w \mid w  \in (A \cup A^{-1})^\ast, w =_G 1 \}
\]
is called the (group-theoretic) \textit{word problem} for $G$, see e.g. \cite{Anisimov1971,Muller1983}. Let $M$ be a monoid with a finite generating set $A$. Translating the above definition of the word problem directly to $M$ does not, in general, yield much insight into the structure of $M$. Duncan \& Gilman \cite[p.\ 522]{Duncan2004}, realising this, introduced a different language. The (monoid) \textit{word problem of $M$ with respect to $A$} is defined as the language
\[
\WP_A^M := \{ u \# v^\trev \mid u, v \in A^\ast, u =_M v\},
\]
where $\#$ is some fixed symbol not in $A$. For a class of languages $\cc$, we say that $M$ has $\cc$-word problem if $\WP_A^M$ is in $\cc$. If $\cc$ is closed under inverse homomorphism, then $M$ having $\cc$-word problem does not depend on the finite generating set $A$ chosen for $M$ \cite[Theorem~5.2]{Duncan2004}. Furthermore, a group has group-theoretic word problem in $\cc$ if and only if $\WP_A^M$ is in $\cc$ \cite[Theorem~3]{Duncan2004}. 

\subsection{Monadic ancestry and super-$\AFL$s}

Many of our theorems will be stated for a special type of classes of languages. Such classes of languages are called \textit{super-$\AFL$s}, and were introduced by Greibach \cite{Greibach1970}.

We follow Book, Jantzen \& Wrathall \cite{Book1982b} in the following definitions. Let $A$ be an alphabet. For each $a \in A$, let $\sigma(a)$ be a language (over any finite alphabet); for every $x, y \in A^\ast$ let $\sigma(xy) = \sigma(x)\sigma(y)$; and for every $L \subseteq A^\ast$, let $\sigma(L) = \bigcup_{w\in L} \sigma(w)$. We then say that $\sigma$ is a \textit{substitution}. For a class $\cc$ of languages, if for every $a \in A$ we have $\sigma(a) \in \cc$, then we say that $\sigma$ is a $\cc$-\textit{substitution}. Let $A$ be an alphabet, and $\sigma$ a substitution on $A$. For every $a \in A$, let $A_a$ denote the smallest finite alphabet such that $\sigma(a) \subseteq A_a^\ast$. Extend $\sigma$ to $A \cup (\bigcup_{a \in A}A_a)$ by defining $\sigma(b) = \{ b \}$ whenever $b \in (\bigcup_{a \in A} A_a) \setminus A$. For $L \subseteq A^\ast$, let $\sigma^1(L) = \sigma(L)$, and let $\sigma^{n+1}(L) = \sigma(\sigma^n(L))$ for $n \geq 1$. Let $\sigma^\infty(L) = \bigcup_{n\geq 0} \sigma^n(L)$. Then we say that $\sigma^\infty$ is an \textit{iterated substitution}. If for every $b \in A \cup (\bigcup_a A_a)$ we have $b \in \sigma(b)$, then we say that $\sigma^\infty$ is a \textit{nested} iterated substitution. We say that $\cc$ is closed under nested iterated substitution if for every $\cc$-substitution $\sigma$ and every $L \in \cc$, we have: if $\sigma^\infty$ is a nested iterated substitution, then $\sigma^\infty(L) \in \cc$. 

Let $\cc$ be a class of languages. We say that $\cc$ is a \textit{super-$\AFL$} if it is an $\AFL$ (i.e. it is closed under homomorphism, inverse homomorphism, intersection with regular languages, union, concatenation, and the Kleene star) and if it closed under nested iterated substitution. As mentioned, the class $\CF$ of context-free languages is a super-$\AFL$ \cite{Kral1970}, \cite[Theorem~2.2]{Book1982b}, as is the class $\IND$ of indexed languages \cite{Aho1968, Engelfriet1985}. On the other hand, the class $\REG$ is clearly not a super-$\AFL$. Indeed, if $\cc$ is a super-$\AFL$, then $\CF \subseteq \cc$, by \cite[Theorem~2.2]{Greibach1970}. Thus $\CF$ is the smallest super-$\AFL$. For more examples and generalisations, we refer the reader to the so-called \textit{hyper-$\AFL$s} defined by Engelfriet \cite{Engelfriet1985}, all of which are super-$\AFL$s. 

A crucial definition, mixing rewriting systems with formal language theory, is the following. Let $\mathcal{C}$ be a class of languages. We say that a rewriting system $\mathcal{R} \subseteq A^\ast \times A^\ast$ is $\mathcal{C}$-\textit{ancestry preserving} if for every $L \in \mathcal{C}$ with $L \subseteq A^\ast$, we have $\langle L \rangle_{\mathcal{R}} \in \mathcal{C}$. A class of languages $\mathcal{C}$ has the \textit{monadic ancestor property} if every monadic $\mathcal{C}$-rewriting system is $\mathcal{C}$-ancestry preserving. It is not hard to see that any $\AFL$ closed under nested iterated substitution has the monadic ancestor property \cite[Theorem~2.2]{Book1982b}). The converse holds, too, though we shall only require the forward direction herein.

\begin{proposition}
Let $\cc$ be an $\AFL$. Then $\cc$ is a super-$\AFL$ if and only if it has the monadic ancestor property. 
\end{proposition}

We refer the reader to Nyberg-Brodda \cite{NybergBrodda2021f} for a full (rather straightforward) proof. Thus we have a combinatorial characterisation of super-$\AFL$s in terms of rewriting systems. We now turn to the main objects of study of this article.

\subsection{Special monoids}

Let $M = \pres{Mon}{A}{w_1 =1, w_2 = 1, \dots, w_i = 1, \dots}$. Then $M$ is called \textit{special}. That is, a monoid is special if it admits a presentation in which the right-hand side of every defining relation is the empty word. Unless stated otherwise, we always consider finitely presented special monoids. The \textit{group of units} $U(M)$ is the subgroup of $M$ consisting of all (two-sided) invertible elements. We say that a word $w \in A^\ast$ is an \textit{invertible} word if it represents an invertible element of the monoid $M$. Arguably the most fundamental result for special monoids is the following, which relates a special monoid to its group of units. 

\begin{theorems}[Makanin, 1966]
Let $M = \pres{Mon}{A}{w_1 =1, w_2 = 1, \dots, w_k = 1}$. Then $U(M)$ is a $k$-relator group. Furthermore, the word and divisibility problems for $M$ reduce to the word problem for $U(M)$.
\end{theorems}

The results were announced in \cite{Makanin1966}, and a proof appeared in Makanin's Ph.D. thesis \cite{Makanin1966b}, see \cite{NybergBrodda2021c} for an English translation by the author of the present article. For the case $k=1$, i.e. special one-relation monoids $\pres{Mon}{A}{w=1}$, and more generally the case when $|w_i| = |w_j|$ for all $i, j \geq 1$, Makanin's result had already been proved by Adian \cite{Adian1960}. In particular, we deduce the decidability of the word problem for any special one-relation monoid $\pres{Mon}{A}{w=1}$, as the word problem for one-relator \textit{groups} $\pres{Gp}{X}{r=1}$ is decidable, a classical result of Magnus' \cite{Magnus1932}. 

We give some definitions and notation which will be used throughout the article. We follow Zhang \cite{Zhang1992} in notation, who produced a significantly shorter rewriting of Makanin's result in terms of rewriting systems. We begin by stating a simple lemma, which appears in most of the arguments used throughout this article, and indeed in the literature in general.\footnote{For example, it appears in Adian \cite[Lemma~3]{Adian1966}, Makanin \cite[Lemma~3]{Makanin1966b}, McNaughton \& Narendran \cite[Lemma~2(2)]{McNaughton1987}, Nivat \cite{Nivat1966}, Lallement \cite{Lallement1974}, Perrin \& Schupp \cite[Lemme]{Perrin1984}, Kobayashi \cite[Corollary~3.3]{Kobayashi2000}, Zhang \cite[Proposition~2.1(4)]{Zhang1992}, and many others.} We shall almost always use it implicitly. 

\begin{lemma}\label{Lem:Fundamental_lemma}
If the words $xy$ and $yz$ are invertible in a monoid $M$, then all three of the words $x, y$, and $z$ are invertible in $M$.
\end{lemma}

An example of a useful application of Lemma~\ref{Lem:Fundamental_lemma} is: if a word $x$ is a prefix of an invertible word $u$, and a suffix of an invertible word $v$, then $x$ is itself invertible.

Consider now an arbitrary finitely presented special monoid
\begin{equation}\label{Eq:Intro_special_monoid}
M = \pres{Mon}{A}{w_1 =1, w_2 = 1, \dots, w_k = 1},
\end{equation}
which shall remain fixed throughout this section. The words $w_1, w_2, \dots, w_k$ will be called the \textit{defining words} of the monoid. We say that an invertible word $u \in A^+$ is \textit{minimal} if none of its non-empty proper prefixes is invertible (in $M$). The set of all minimal words forms a biprefix code, denoted $\mathfrak{M}$. Obviously any invertible word factors (uniquely) as a product of minimal words. 

Every defining word $w_i$ for $1 \leq i \leq k$ is an invertible word, as $w_i =_M 1$. Hence, we can uniquely factor every $w_i$ into minimal words as $w_i \equiv w_{i,1} w_{i,2} \cdots w_{i,\ell_i}$, where $w_{i,j} \in \mathfrak{M}$ for $1 \leq j \leq \ell_i$. The set of all minimal words arising in this way shall be denoted $\Lambda$, and called the set of \textit{presentation pieces} of \eqref{Eq:Intro_special_monoid}. That is, 
\[
\Lambda = \bigcup_{i=1}^{k} \bigcup_{j=1}^{\ell_i} \{ w_{i,j} \} \subseteq A^\ast.
\]

We let $\Delta$ denote the set of all minimal words $\delta \in \mathfrak{M}$ satisfying: there exists some $\lambda \in \Lambda$ with $\delta =_M \lambda$ and $|\delta| \leq |\lambda|$. The set $\Delta$ is called the set of \textit{invertible pieces} of the presentation.\footnote{This definition of $\Delta$ is slightly more restrictive than Zhang's. Namely, Zhang defines $\Delta$ as the set of minimal words $\delta \in A^\ast$ such that there exists some $\lambda \in \Lambda$ with $\delta =_M \lambda$ and $|\delta| \leq \max_{\lambda \in \Lambda} |\lambda|$. However, the only place in \cite{Zhang1992} where the $\max_{\lambda \in \Lambda}$ is used is in the proof of \cite[Prop~2.3]{Zhang1992}, and there one can with no further changes necessary substitute the present definition. This is done in detail in the author's Ph.D. thesis \cite[\S1.3]{Thesis}. Thus, all of Zhang's results can be applied here.} As $\Delta \subseteq \mathfrak{M}$, no elements of $\Delta$ overlap non-trivially. In particular, $\Delta$ is a biprefix code. Furthermore, note that $\Lambda \subseteq \Delta$ and that $\langle \Delta \rangle = \langle \Lambda \rangle$, as submonoids of $M$. We partition $\Delta$ according to which words in $\Delta$ represent the same elements of $M$. This partition of $\Delta$, i.e. the partition of $\Delta$ induced by the equivalence relation $=_M$, will be denoted $\Delta_1 \cup \Delta_2 \cup \cdots \cup \Delta_\nu$. Let $X = \{ x_1, \dots, x_\nu \}$ be a set of new symbols, and let $\phi \colon \Delta^\ast \to X^\ast$ be the map induced by $\delta \mapsto x_i$ when $\delta \in \Delta_i$. This is a well-defined homomorphism, as $\Delta$ is a biprefix code. One can show (see \cite[Theorem~3.7]{Zhang1992}) that 
\begin{equation}\label{Eq:Pres_for_UM}
\pres{Gp}{X}{ \{ \phi(w_i) = 1 \: (1 \leq i \leq k) \} }
\end{equation}
is a ($k$-relator) group presentation for the group of units $U(M)$ of $M$. In particular, we have $\langle \Delta \rangle = \langle \Lambda \rangle = U(M)$. In fact, for $u, v \in \Delta^\ast$ we have $u =_M v$ if and only if $\phi(u) =_{U(M)} \phi(v)$, see \cite[Lemmas~3.1 and 3.6]{Zhang1992}. We will find it convenient to consider $\Delta$ as a finite generating set for $U(M)$, and consider the language
\begin{equation}\label{Eq: WPDeltaUM}
\WP_\Delta^{U(M)} = \{ u \# v^\trev \mid u, v \in \Delta^\ast, \phi(u) =_{U(M)} \phi(v)\} = \WP_A^M \cap \Delta^\ast \# (\Delta^\trev)^\ast.
\end{equation}

Zhang \cite{Zhang1992} introduced a rewriting system for studying special monoids, which we shall also find use for in this article. For any total order $\prec$ on $A$, let $<_{s}$ denote the short-lex order on $A$ induced by $\prec$. Define the rewriting system $S = S(M)$ as:
\begin{equation}\label{Eq:S}
S = \{ (u \to v) \mid u, v \in \Delta^\ast \colon u =_M v \textnormal{ and } u >_s v\}.
\end{equation}
Zhang \cite[Proposition~3.2]{Zhang1992} proves that this (generally infinite) system defines $M$, and is furthermore a complete rewriting system. 

Now, in general, neither the set $\Delta$ nor the presentation \eqref{Eq:Pres_for_UM} are effectively computable from \eqref{Eq:Intro_special_monoid}. However, if $k=1$, then there is an algorithm, called \textit{Adian's (overlap) algorithm}, which computes both $\Delta$ and \eqref{Eq:Pres_for_UM}. In fact, in this case, we also have $\Lambda = \Delta$ and $\nu=1$, which follows from Magnus' \textit{Freiheitssatz}. We refer the reader to \cite{Lallement1974, Zhang1992b, Gray2021} for details. Adian's algorithm can fail to produce the correct output already when $k=2$. Indeed, in general the problem of computing \eqref{Eq:Pres_for_UM} is undecidable. A procedure for when $k\geq 1$, implicit in Makanin's Ph.D. thesis, is described in the author's Ph.D. thesis \cite{Thesis}. This does not always terminate, but when it does it always outputs \eqref{Eq:Pres_for_UM} and $\Delta$ correctly.

We make one final definition. If $\lambda \in \Lambda$ is a presentation piece, then we say that $\lambda$ is obtained from itself by the \textit{piece-generating operation}, and inductively, we say:
 
\begin{itemize}[label=($*$)]
\item Suppose that $w \equiv h_1 \delta_{i,1} \delta_{i,2} \cdots \delta_{i,p} h_2$ is obtained from $\lambda$ by the piece-generating operation, where $p \geq 0$ and $h_1, h_2$ are non-empty, and $\delta_{i,j} \in \Delta$ for every such $\delta_{i,j}$. Suppose then that $w' \equiv h_1 \delta_{j,1} \delta_{j,2} \cdots \delta_{j,t} h_2$, with $t \geq 0$, that $|w'| \leq |w|$, and that $\delta_{i,1} \delta_{i,2} \cdots \delta_{i,p} =_M \delta_{j,1} \delta_{j,2} \cdots \delta_{j,t}$. Then $w'$ is also said to be obtained from $\lambda$ by the piece-generating operation.
\end{itemize} 

If $w \in A^\ast$ can be obtained from $\lambda$ by the piece-generating operation, then we denote this by $w \in [\lambda]^\downarrow$. It is easy to see that for any $w \in [\lambda]^\downarrow$, we have that $w$ is a minimal word, i.e. $w \in \mathfrak{M}$ (see e.g. the second half of the proof of Lemma~\ref{Lem:dD_is_invertibles}). We have $\lambda \in [\lambda]^\downarrow$, and $|w| \leq |\lambda|$ for every $w \in [\lambda]^\downarrow$, so in particular for every $\lambda \in \Lambda$ the set $[\lambda]^\downarrow$ is finite. In general, for $\lambda_1, \lambda_2 \in \Lambda$ we can have $[\lambda_1]^\downarrow \cap [\lambda_2]^\downarrow \neq \varnothing$ even when $\lambda_1 \not\equiv \lambda_2$. We remark the following useful fact: if $\lambda \in \Lambda$ and $\delta \in \Delta$ are such that $\lambda \equiv h_1 w h_2$ and $\delta \equiv h_1 w' h_2$, with $h_1, h_2 \in A^+$ and $w \xra{S} w'$, then $\delta \in [\lambda]^\downarrow$. The converse does not, in general, hold.

\begin{example}
Let $M_1 = \pres{Mon}{a,b,c}{ab^3c = 1, b^2 = 1}$. Then it is not hard to show (e.g. by using a finite complete rewriting system for $M_1$) that $\Lambda = \{ b, ab^3c\}$, while $\Delta = \{ b, abc, ab^3c\}$. Then $abc \in \Delta$ is obtained from the piece $ab^3c \in \Lambda$ by the generating operation, as $b^3 =_{M_1} b$ and $|abc| \leq |ab^3c|$. That is, $abc \in [ab^3c]^\downarrow$.
\end{example}

Otto~\&~Zhang \cite[Theorem~5.2]{Otto1991} proved the following ``normal form lemma''.

\begin{lemma}[Otto \& Zhang]\label{Lem:Zhangs_lemma}
Let $M$ be as given in \eqref{Eq:Intro_special_monoid}, and let $u, v \in A^\ast$ be such that $u =_M v$. Then we can uniquely factorise $u$ and $v$ as
\begin{align*}
u \equiv u_0 a_1 u_1 \cdots a_m u_m, \qquad \textnormal{ and} \qquad v \equiv v_0 a_1 v_1 \cdots a_m v_m,
\end{align*}
respectively, where for all $0 \leq i \leq m$ we have $a_i \in A$ and 
\begin{enumerate}[label=(\arabic*)]
\item $u_i =_M v_i$;
\item $u_i$ is a maximal invertible factor of $u$.
\item $v_i$ is a maximal invertible factor of $v$.
\end{enumerate}
\end{lemma}

Here, a maximal invertible subword of $w$ is one which is not properly contained in any other invertible subword of $w$. Lemma~\ref{Lem:Zhangs_lemma} tells us, roughly speaking, the importance of understand the equality of invertible words for understanding $\WP_A^M$. For this reason, we will begin by studying the invertible words of special monoids.

\section{Invertible elements}

\noindent Throughout this section, if not explicitly stated otherwise, we let
\begin{equation}\label{Eq:Special_monoid_pres}
M = \pres{Mon}{A}{w_1 = 1, w_2 = 1, \dots, w_k =1}
\end{equation}
be a fixed special monoid. We assume $|A|<\infty$. Let $\Delta$ be the set of minimal invertible pieces of $M$, and let $\phi \colon \Delta^\ast \to X^\ast$ be the associated homomorphism. In general, if $w \in A^\ast$ is invertible, it need not be the case that $w \in \Delta^\ast$. For example, in the bicyclic monoid $B = \pres{Mon}{b,c}{bc=1}$, we have $\Delta = \{ bc\}$, but $b^n c^n \not\in \Delta^\ast$ is invertible for every $n \geq 0$. We remark as an aside that it is easy to prove that every non-empty invertible word contains a piece as a subword. The aim of this section is to give an explicit description of the invertible words in terms of $\Delta$, in order to study their language-theoretic properties. 

\subsection{Generalised pieces} Let $w \in A^\ast$ be a minimal word such that there exists some $\delta \in \Delta$ with $w \xra{S} \delta$. Then we say that $w$ is a \textit{generalised piece}. Recall that $S = S(M)$ is the system \eqref{Eq:S}. The collection of all generalised pieces is denoted $\dD$. Note that $\Delta \subseteq \dD$. Furthermore, as $\dD$ consists of minimal words, we conclude that $\dD$ is a biprefix code as a subset of $A^\ast$. The following elementary properties are easy to prove, with the same proof, \textit{mutatis mutandis}, as \cite[Proposition~2.1]{Zhang1992}.

\begin{lemma}\label{Lem:Basic_properties}
Let $x, y, z \in A^\ast$. Then (1) $xy, x \in \dD^\ast$ implies $y \in \dD^\ast$. (2) $yz, z \in \dD^\ast$ implies $y \in \dD^\ast$. (3) Suppose that $xy \in \dD^\ast$. If either $x$ or $y$ is invertible, then $x \in \dD^\ast$ and $y \in \dD^\ast$. (4) Suppose that $xy \in \dD^\ast$ and $yz \in \dD^\ast$. Then $x, y, z \in \dD^\ast$. 
\end{lemma}
\begin{proof}
Statements (1) and (2) follow directly from the fact that $\dD$ is a biprefix code. For (3), suppose $x$ is invertible, and let $xy \equiv x_1 x_2 \cdots x_m$ with $x_i \in \dD$ where $1 \leq i \leq m$. Then $x \equiv x_1 x_2 \cdots x_{\ell-1} E$, where $E$ is a prefix of $x_\ell$ for some $\ell \leq m$. Write $x_\ell \equiv EF$ where $F \in A^\ast$. As $x$ and $x_\ell$ are invertible, it follows that $E$ is invertible by overlaps. As $x_\ell$ is minimal, we must thus have that either $E$ is empty, or else $E$ is all of $x_\ell$. In either case, $x \in \dD^\ast$. By (1), we conclude $y \in \dD^\ast$. Symmetrically, the results hold when $y$ is invertible. For (4), as we have $xy \in \dD^\ast$ and $yz \in \dD^\ast$, we have that $xy$ and $yz$ are invertible. Hence $y$ is invertible. By (3), we have $x, y, z \in \dD^\ast$. 
\end{proof}

Lemma~\ref{Lem:Basic_properties}(4) ensures that the types of overlap arguments possible for words over $\Delta^\ast$ are also possible for $\dD^\ast$. We shall use this implicitly throughout the remainder of the article. Now, $\dD$ is a set whose elements are, in a certain sense, controlled by elements of $\Delta$. The main benefit of introducing this set is the following lemma.

\begin{lemma}\label{Lem:dD_is_invertibles}
A word $w \in A^\ast$ is invertible if and only if $w \in \dD^\ast$. 
\end{lemma}
\begin{proof}
Any element of $\dD^\ast$ is clearly invertible. For the converse, assume $w \in A^\ast$ is invertible. By \cite[Lemma~3.4]{Zhang1992} there exists some least $n \geq 0$ and a $D \in \Delta^\ast$ such that $w \xra{S}^n D$. We will prove the claim by induction on $n$. The base case $n=0$ is clear, for then $w \equiv D \in \Delta^\ast \subseteq \dD^\ast$. Assume that the claim is true for some $n \geq 0$, and let $w$ be such that $w \xra{S}^{n+1} D$. Then there exists some $w_1 \in A^\ast$ such that $w \xr{S} w_1$ and $w_1 \xra{S}^n D$. As $w =_M w_1$, the word $w_1$ is invertible and by the inductive hypothesis $w_1 \in \dD^\ast$. Write $w_1 \equiv \overline{\delta}_0 \overline{\delta}_1 \cdots \overline{\delta}_k$ where $\overline{\delta}_i \in \dD$ for $0 \leq i \leq k$. As $w \xr{S} w_1$, there exists some $(\ell, r) \in S$ and words $u, v \in A^\ast$ such that $w \equiv u \ell v$ and $w_1 \equiv urv$. We subdivide into two cases, depending on whether $r$ contains as a subword one of the $\overline{\delta}_i$ or not. 

In the first case, we assume the fixed subword $r$ of $w_1$ contains some $\overline{\delta}_i$, where $0 \leq i \leq k$, as a subword. Then, as $\dD$ is a biprefix code and $r \in \Delta^\ast \subseteq \dD^\ast$, we must have that $w_1 \equiv ErF$, where $E, F \in \dD^\ast$. Hence $w \equiv E\ell F$, and as $\ell \in \Delta^\ast$, we have $w \in \dD^\ast$. In the second case, this fixed subword $r$ does not contain any $\overline{\delta}_i$ as a subword. We deal with two separate subcases, depending on whether $r$ is empty or not. 

First, if $r \equiv \varepsilon$, then there exists $0 \leq i \leq k$ we can write $u \equiv \overline{\delta}_0 \cdots \overline{\delta}_{i-1}\overline{\delta}'_i$ and $v \equiv \overline{\delta}''_i \overline{\delta}_{i+1} \cdots \overline{\delta}_k$, where $\overline{\delta}'_i, \overline{\delta}''_i \in A^\ast$ are such that $\overline{\delta}'_i \overline{\delta}''_i \equiv \overline{\delta}_i$; and such that 
\begin{equation}\label{Eq:w=ulv=delta...}
w \equiv u\ell v \equiv \overline{\delta}_0 \cdots \overline{\delta}_{i-1}(\overline{\delta}_i'\ell \overline{\delta}''_i) \overline{\delta}_{i+1} \cdots \overline{\delta}_k.
\end{equation}
Assume $|\overline{\delta}'_i|\cdot |\overline{\delta}''_i| = 0$, i.e. at least one of $\overline{\delta}'_i, \overline{\delta}''_i$ is empty. If $\overline{\delta}'_i \equiv \varepsilon$, then $\overline{\delta}''_i \equiv \overline{\delta}_i \in \dD$, and as $\ell \in \dD^\ast$, we have $w \in \dD^\ast$ by \eqref{Eq:w=ulv=delta...}. The case $\overline{\delta}''_i \equiv \varepsilon$ is entirely symmetrical. Thus assume instead that $|\overline{\delta}'_i|\cdot |\overline{\delta}''_i| > 0$. We claim that no non-empty prefix of $\overline{\delta}'_i \ell \overline{\delta}''_i$ is invertible, i.e. this word is minimal. By minimality of $\overline{\delta}'_i\overline{\delta}''_i \in \dD$, no proper non-empty prefix or suffix of this word is invertible; thus if some prefix of $\overline{\delta}'_i \ell \overline{\delta}''_i$ were invertible, then it is of the form $\overline{\delta}'_i \ell_1$, where $\ell_1 \in A^+$ is some non-empty proper prefix of $\ell$. Thus $\ell_1$ is left invertible, being a suffix of the invertible $\overline{\delta}'_i \ell_1$, but also right invertible, being a prefix of $\ell$. It follows that $\ell_1$ is invertible. As $\overline{\delta}_i' \ell_1$ is invertible, thus $\overline{\delta}_i'$ is invertible, which contradicts the minimality of $\overline{\delta}_i \equiv \overline{\delta}'_i \overline{\delta}''_i$ as $|\overline{\delta}''_i| > 0$. It follows that $\overline{\delta}'_i \ell \overline{\delta}''_i$ is minimal. As $\overline{\delta}'_i \ell \overline{\delta}''_i$ it is clearly invertible, being equal in $M$ to $\overline{\delta}'_i \overline{\delta}''_i \in \dD$ by virtue of $\ell =_M 1$, we have that $\overline{\delta}'_i \ell \overline{\delta}''_i \in \dD$. By \eqref{Eq:w=ulv=delta...}, we have $w \in \dD^\ast$. 

The case $r \not\equiv \varepsilon$ uses very similar reasoning. We omit the proof for brevity. 
\end{proof}

\begin{proposition}
$\dD = \mathfrak{M}$. 
\end{proposition}
\begin{proof}
Clearly, $\mathfrak{M}^\ast$ is the set of invertible words of $M$, so $\mathfrak{M}^\ast = \dD^\ast$ by Lemma~\ref{Lem:dD_is_invertibles}. As $\mathfrak{M}$ and $\dD$ are both biprefix codes, we thus necessarily have $\dD = \mathfrak{M}$. 
\end{proof}

The description of $\mathfrak{M}$ as $\dD$ gives us access to an explicit description of the minimal words in terms of the objects $\Delta$ and $\xra{S}$, which will be useful. For some more properties of the set $\dD$, we refer the reader to Chapter~3 of the author's Ph.D. thesis \cite{Thesis}.

\subsection{Controlling the pieces}\label{Subsec:Controlling_pres} We will now demonstrate certain manipulations of the pieces $\Delta$ of a presentation, in order to gain sufficient control over the set $\dD$. One desirable property, which would simplify a bulk of reasoning, would be the non-existence of elements of $\Delta$ appearing as proper subwords of other pieces in $\Delta$, i.e. that $\Delta$ is an \textit{infix} code.\footnote{Garreta \& Gray \cite{Garreta2019} call this condition (C1), but do not directly study it.} We do not know if every special monoid admits a presentation where $\Delta$ is an infix code (see Question~\ref{Quest:all_infix_free?}). Instead, we introduce and consider a weaker property, which shall be sufficient for our purposes. Let $\delta \in \Delta$ be a piece. If $\delta \equiv h_1 w h_2$ for some $h_1, h_2 \in A^+$ and $w \in \Delta^+$, then we say that $w$ is a \textit{subpiece} (of $\delta$). If $w$ is a subpiece and $|w|=1$, then we say that $w$ is a \textit{small} subpiece. A subpiece which is not small is called \textit{large}. The special monoid presentation \eqref{Eq:Special_monoid_pres} is said to satisfy the \textit{small subpiece condition} if all subpieces of pieces in $\Delta$ are small. 

\begin{example}
If $M_2 = \pres{Mon}{a,b,c}{ab^2c =1, b = 1}$, then one can show that the pieces are $\Delta = \{ ab^2c, abc, ac, b \}$. Thus there are three subpieces; $b$ appearing as a subword of $abc$ and of $ab^2c$, and $b^2$ appearing as a subword of $ab^2c$. The first two subpieces are small, but as $|b^2|=2$, this last subpiece is not small. 
\end{example}
We remark that, in the above two examples, we obviously have $M_1 \cong M_2$. Hence, the property of satisfying the small subpiece condition is quite strongly tied to the presentation chosen for the monoid. In this section we will prove the following:

\begin{proposition}\label{Prop:Every_admits_small_subpiece}
Every finitely presented special monoid admits a presentation satisfying the small subpiece condition. 
\end{proposition}

Before proving this, we state a useful lemma, proved by Makanin \cite[Lemma~12]{Makanin1966b}. 

\begin{lemma}[Makanin's Lemma]\label{Lem:Makanins_lemma}
Let $M$ be a special monoid given by the finite presentation
\[
M= \pres{Mon}{A}{w_1 = 1, \dots, w_k = 1}.
\]
with presentation pieces $\Lambda$. Fix some $1 \leq i \leq k$, and let $w_i \equiv \lambda_1 \lambda_2 \cdots \lambda_\ell$ with $\lambda_j \in \Lambda$ for $1 \leq j \leq \ell$. Suppose that $\delta \in \Delta$ is such that $\delta \in [\lambda_p]^\downarrow$ for some fixed $1 \leq p \leq \ell$. Let $w_i' \equiv \lambda_1 \lambda_2 \cdots \lambda_{p-1} \delta \lambda_{p+1} \cdots \lambda_\ell$, and let 
\[
M' = \pres{Mon}{A}{w_{1}=1, \dots, w_i' = 1, \dots, w_k = 1}.
\]
Then $M \cong M'$ via the identity map, i.e. the presentations define the same congruence on $A^\ast$. Furthermore, the factorisation in $M'$ of the defining word $w'_i$ into minimal invertible factors is obtained by replacing $\lambda_p$ with $\delta$ in the factors of the factorisation of $w_i$, and the factorisation of the defining word $w_j$ ($j \neq i$, $1 \leq j \leq k$) is identical to its factorisation in $M$.
\end{lemma}

For a full discussion of the validity of the translation of the lemma into the language of minimal invertible pieces $\Delta$ (rather than Makanin's original ``$c$-words''), see \S1.3 of the author's Ph.D. thesis \cite{Thesis}. We will now give an example for how Makanin's lemma can be applied to remove large subpieces. The idea in the example is the same as the general idea which will be used in the proof of Proposition~\ref{Prop:Every_admits_small_subpiece}. 

\begin{example}
Let $M_3 = \pres{Mon}{a,b}{abaabbab = 1}$. By Adian's algorithm (see Lallement \cite{Lallement1974}), the defining word factors into invertible pieces as $(ab)(aabb)(ab)$, so $\Delta = \Lambda = \{ ab, aabb \}$. Thus $ab \in \Delta^+$ is a large subpiece of $aabb \in \Delta$. We will replace this large subpiece $ab$ by a small subpiece $p$.

Let $p$ be any new symbol, and introduce the defining relation $p = ab$ to the presentation via a Tietze transformation, giving $\pres{Mon}{a,b,p}{abaabbab = 1, p=ab}$. It is clear that $aabb \cdot ab$ is an inverse of $ab$, so from $p=_{M_3} ab$ it thus follows that $M_3$ is isomorphic to 
\[
\pres{Mon}{a,b,p}{abaabbab = 1, p=ab, p(aabb \cdot ab) = 1, (aabb \cdot ab)p = 1}.
\]
As the fact that both $p$ and $ab$ are inverses of $aabbab$ follows from the two added relations, the relation $p=ab$ follows from these relations; thus we can remove $p=ab$, and find that $M_3$ is isomorphic to the special monoid
\begin{equation}\label{Eq:ab(aabb)(ab)=1}
\pres{Mon}{a,b}{(ab)(aabb)(ab) = 1, p(aabb)(ab)=1, (aabb)(ab)p = 1}.
\end{equation}
It is clear that each of $p, aabb, apb$, and $ab$ is a minimal invertible piece of this presentation. As $p=_{M_3} ab$, we have that the piece $apc$ is obtained from $aabb$ by the piece-generating operation. Thus, by Makanin's Lemma, we may replace $aabc$ by $apc$ in \eqref{Eq:ab(aabb)(ab)=1} without changing the monoid defined by it. Thus
\[
M_3 \cong \pres{Mon}{a,b,p}{(ab)(apb)(ab)=1, p(apb)(ab) = 1, (apb)(ab)p = 1}.
\]
It is now not difficult to show that, for this new presentation, $\Delta = \Lambda = \{ p, ab, apb\}$. Thus, this new presentation satisfies the small subpiece condition. 
\end{example}

We now generalise the above example to the general case. First, for any set $S \subseteq A^\ast$, we let $\omega(S)$ denote the natural number $\sum_{s \in S}(|s|-1)$. If $S = \Lambda$, where $\Lambda$ is the set of presentation pieces of \eqref{Eq:Special_monoid_pres}, then in a loose sense $\omega(\Lambda)$ is a measure of the ``complexity'' of the pieces of the presentation -- we remark that $M$ is right cancellative if and only if $\omega(\Lambda)=0$ (see \S\ref{Subsec:Other_classes}). We shall, in the subsequent proof of Proposition~\ref{Prop:Every_admits_small_subpiece}, use several operations on the presentation \eqref{Eq:Special_monoid_pres}, each of which reduces or does not increase $\omega(\Lambda)$. In particular, if a presentation \eqref{Eq:Special_monoid_pres}, with presentation pieces $\Lambda$, does not satisfy the small subpiece condition, we will show that we find a new presentation with presentation pieces $\Lambda'$ such that $\omega(\Lambda')<\omega(\Lambda)$. The proof will then be complete by induction. 

Before we can realise the above idea in practice, we remark that it is not difficult to construct special monoids in which no presentation piece has a large subpiece, but there is some piece with a large subpiece. We begin with a lemma to remedy this, by showing that we can always find a presentation where large subpieces are ``brought to light'' in the presentation pieces.

\begin{lemma}\label{Lem:Admits_either_small_or_presentationpiece}
Let $M$ be as in \eqref{Eq:Special_monoid_pres}, with presentation pieces $\Lambda$. Then $M$ admits a special monoid presentation, with presentation pieces $\Lambda'$, such that either this presentation satisfies the small subpiece condition; or else there is a presentation piece $\lambda \in \Lambda'$ containing a large subpiece, and $\omega(\Lambda') \leq \omega(\Lambda)$.
\end{lemma}
\begin{proof}
If \eqref{Eq:Special_monoid_pres} satisfies the small subpiece condition, then we are done, so suppose that it does not. Let $\delta \in \Delta$ and $w \in \Delta^+$ be such that $w$ is a large subpiece of $\delta$. If $\delta \in \Lambda$, then we are done. If $\delta \not\in \Lambda$, then there is some $\lambda_0 \in \Lambda$ such that $\delta =_M \lambda_0$. Fix such a $\lambda_0$. Then there are words $u_0, u_1, \dots, u_n \in A^\ast$ and a sequence 
\begin{equation}\label{Eq:delta-to-lambda}
\delta \equiv u_0 \lr{M} u_1 \lr{M} \cdots \lr{M} u_{n-1} \lr{M} u_n \equiv \lambda_0.
\end{equation}
In the rewriting \eqref{Eq:delta-to-lambda}, suppose (without loss of generality, by symmetry) that the first letter of $\delta$ is affected (in the sense of Novikov \cite[I.\S1]{Novikov1955} and Adian \cite{Adian1966}) before the last letter is. Suppose the first time, if any, this happens is in the rewriting $u_i \xr{M} u_{i+1}$. Then $u_i \equiv v_0 \delta' v_1$, where $v_0 =_M v_1 =_M 1$ and $\delta' =_M \delta$. The rewriting $u_i \xr{M} u_{i+1}$ affects the first letter of $\delta'$ (which is the same as the first letter of $\delta$) by deleting a defining relation $w_j$ ($1 \leq j \leq k$), and therefore must, by minimality of $\delta'$, and invertibility of $v_0, v_1$, be such that $\delta'$ is one of the minimal invertible pieces in the factorisation of this $w_j$; thus $\delta' \in \Lambda$.

We conclude that for our chosen $\delta$, we can find a piece $\lambda \in \Lambda$ such that $\delta =_M \lambda$, and there is a rewriting $\delta \lra{M} \lambda$ which does not affect the first or last letter of $\delta$. Indeed, we can take $\lambda \equiv \delta'$ as above if the first letter of $\delta$ is affected in \eqref{Eq:delta-to-lambda}; otherwise, we can take $\lambda \equiv \lambda_0$. In either case, pick the longest $\lambda$ with the given property. Then there exist $h_1, h_2 \in A^+$ with $\delta \equiv h_1 w h_2$ and $\lambda \equiv h_1 w' h_2$ with $w =_M w'$. Thus there is some $W \in A^\ast$ with $w \xra{S} W$ and $w' \xra{S} W$. Hence $|w| \geq |W|$ and $|w'| \geq |W|$. As $\lambda$ was chosen longest, we also have $|w| \leq |w'|$ (for otherwise $\delta \not\in \Delta$). 

Now, if $|w'| = |W|$, then also $|w| = |W|$. Thus the sequences of rules 
\begin{align*}
(s_{1,1} , s_{1,2}), (s_{2,1} , s_{2,2}), \dots, (s_{m,1} , s_{m,2}) &\in S \\
(s'_{1,1} , s'_{1,2}), (s'_{2,1} , s'_{2,2}), \dots, (s'_{\ell,1} , s'_{\ell,2}) &\in S
\end{align*}
transforming $w \xra{S} W$ resp. $w' \xra{S} W$ satisfies $|s_{i,1}| = |s_{i,2}|$ resp. $|s'_{j,1}| = |s'_{j,2}|$ for all $1 \leq i \leq m$ resp. $1 \leq j \leq \ell$. Thus, by composing the sequence of rules rewriting $w'$ to $W$ with the reverse of the sequence of rules rewriting $w$ to $W$, we find a sequence of applications of the piece-generating operation rewriting $h_1 w' h_2$ to $h_1 w h_2$. In other words, $h_1 w h_2 \in [h_1 w' h_2]^\downarrow$, i.e. $\delta \in [\lambda]^\downarrow$. By Makanin's Lemma, we may everywhere replace $\lambda$ by $\delta$ without changing the presentation; in the resulting presentation, whose presentation pieces will be denoted $\Lambda'$, we have $\delta \in \Lambda'$, and $\delta$ contains a large subpiece. As $|\delta| = |\lambda'|$, we have $\omega(\Lambda') = \omega(\Lambda)$, and we are done.

Suppose instead that $|w'|>|W|$. We have $\delta' :\equiv h_1 W h_2 \in [\lambda]^\downarrow$. By Makanin's Lemma, we may everywhere replace $\lambda$ with $\delta'$ in \eqref{Eq:Special_monoid_pres} without changing the monoid $M$. Let $\Lambda'$ be the new presentation pieces of this presentation. Then, as $\lambda$ was chosen longest and $|\delta'| < |\lambda|$, we have $\omega(\Lambda') < \omega(\Lambda)$, as $\Lambda' = (\Lambda - \{\lambda\}) \cup \{ \delta' \}$ by the second part of Makanin's Lemma. We may thus repeat the above proof for the new presentation, and are done by induction.
\end{proof}

Having proved this rather technical lemma, we can proceed with our main proof.

\begin{proof}[Proof of Proposition~\ref{Prop:Every_admits_small_subpiece}]
Suppose $M$ is given by a presentation \eqref{Eq:A_pres_satisfying_small_subpiece}, with presentation pieces $\Lambda_0$. Then $M$ admits a presentation 
\begin{equation}\label{Eq:Pres_that_is_nice_1}
\pres{Mon}{A}{w_1 = 1, w_2 =1, \dots, w_k = 1}
\end{equation}
satisfying the conclusions of Lemma~\ref{Lem:Admits_either_small_or_presentationpiece}, with presentation pieces $\Lambda$ resp. pieces $\Delta$, and such that $\omega(\Lambda) \leq \omega(\Lambda_0)$. If \eqref{Eq:Pres_that_is_nice_1} satisfies the small subpiece condition, we are done, so assume the second part of Lemma~\ref{Lem:Admits_either_small_or_presentationpiece} holds, and let $\lambda \in \Lambda$ be a presentation piece such that $w \in \Delta^+$ is a large subpiece of $\lambda$. Write $\lambda \equiv h_1 w h_2$ with $h_1, h_2 \in A^+$. Introduce a new symbol $p$, disjoint from $A$, and add by way of Tietze transformation the relation $p=w$ to the presentation \eqref{Eq:Pres_that_is_nice_1}. The resulting presentation is not special. However, as $w$ is invertible, there exists some $w' \in \Lambda^+$ such that $ww' =_M w'w =_M 1$. Hence also $pw' =_M w'p =_M 1$. We add these relations to the presentation. As inverses in a group are unique, and as $p$ and $w$ are both invertible words, we find the relation $p=w$ redundant. We remove it by a Tietze transformation, resulting in a new special presentation:
\begin{equation}\label{Eq:Pres_that_is_nice_2}
M' = \pres{Mon}{A \cup \{ p \}}{w_1 = 1, \dots, w_k = 1, pw' = 1, w'p = 1}.
\end{equation}
Now the map induced by $a \mapsto a$ for all $a \in A$ extends to an isomorphism from $M$ to $M'$. Thus the factorisation of $w'$ and the $w_i$, for $1 \leq i \leq k$, into minimal invertible pieces is the same in $M'$ as in $M$. Clearly, $p$ is invertible. It follows that the set $\Lambda'$ of presentation pieces of \eqref{Eq:Pres_that_is_nice_2} is precisely $\Lambda' = \Lambda \cup \{ p \}$.

From the presentation piece $\lambda \equiv h_1 w h_2 \in \Lambda \subset \Lambda'$ in $M'$ we can by the piece-generating operation obtain the piece $\delta := h_1 p h_2$, as $p, w' \in (\Lambda')^\ast$ satisfy $p =_{M'} w$ and $|p| < |w|$. That is, $\delta \in [\lambda]^\downarrow$. By Makanin's Lemma, we can thus in the factorisations of the defining words in \eqref{Eq:Pres_that_is_nice_2} replace $\lambda$ by without changing the monoid. Let $w_i'$ denote the word obtained by this replacement from $w_i$ (for $1 \leq i \leq k$), and $w''$ the word from $w'$. We find a new presentation 
\begin{equation}\label{Eq:Pres_that_is_nice_3}
M'' = \pres{Mon}{A \cup \{ p \}}{w_1' = 1, \dots, w_k' = 1, pw'' =1, w''p = 1}.
\end{equation}
Let $\Lambda''$ denote the presentation pieces of \eqref{Eq:Pres_that_is_nice_3}. As $|p| = 1$ and $|w|>1$, it follows that $\delta := h_1 p h_2$ satisfies $|\delta|<|\lambda|$. By the second part of Makanin's Lemma, $\delta \in \Lambda''$, and the other presentation pieces of \eqref{Eq:Pres_that_is_nice_3} are presentation pieces of \eqref{Eq:Pres_that_is_nice_2}, i.e. in $\Lambda'$. Thus $\omega(\Lambda'') < \omega(\Lambda')$. In particular, we find 
\[
\omega(\Lambda'') < \omega(\Lambda') = \sum_{\lambda' \in \Lambda \cup \{ p \}} (|\lambda'|-1) = (1-1) + \sum_{\lambda' \in \Lambda}(|\lambda'|-1) = \omega(\Lambda).
\]
Thus, repeating the above proof starting with the presentation \eqref{Eq:Pres_that_is_nice_3}, either \eqref{Eq:Pres_that_is_nice_3} satisfies the small subpiece condition, or else we obtain a presentation $M'''$ with presentation pieces $\Lambda'''$ satisfying $\omega(\Lambda''') < \omega(\Lambda'')$, etc. We conclude by induction on $\omega$ that there is some $n \geq 0$ and a presentation $M^{(n)}$ with pieces $\Delta^{(n)}$ such that no piece $\delta \in \Delta^{(n)}$ has a large subpiece; that is, $M^{(n)}$ satisfies the small subpiece condition, and defines $M$.
\end{proof}

Presentations satisfying the small subpiece condition will now prove crucial for describing the language of representatives of pieces in special monoids. 

\subsection{Representatives of pieces} Recall that $X$ is a set in bijective correspondence of cardinality the size $\nu$ of the partition of $\Delta$ as $\Delta_1 \cup \Delta_2 \cup \cdots \cup \Delta_\nu$ into pieces equal to each other in $M$, and that $\phi \colon \Delta^\ast \to X^\ast$ is the canonical surjective homomorphism. For a word $w \in \Delta^\ast$, we define the language of \textit{$\Delta$-representatives} of $w$ as the set 
\begin{equation}\label{Eq:Drep_def}
\DRep_A^M(w) := \{ v \in \Delta^\ast \mid v =_M w\} = \phi^{-1}\left( \Rep_X^{U(M)}(\phi(w)) \right).
\end{equation}
From the rightmost representation in \eqref{Eq:Drep_def}, it follows that if $\cc$ is a class of languages closed under inverse homomorphisms, then $\DRep_A^M(w) \in \cc$ if and only if $\Rep_X^{U(M)}(\phi(w)) \in \cc$. For example, in the bicyclic monoid $B = \pres{Mon}{b,c}{bc=1}$, with $\Delta = \{ bc\}$, $X = \{ x_1 \}$, and $U(B) = \pres{Mon}{x_1}{x_1 = 1}$, we have 
\[
\DRep_{\{b, c\}}^B(1) = \phi^{-1}\left( \Rep_{X}^{U(B)}(\phi(1))\right) = \phi^{-1}\left( x_1^\ast \right) = (bc)^\ast.
\]
Thus $\DRep_{\{b, c \}}^B(1)$ is regular, as $U(B)$ is a group with regular word problem.

The idea of the following section is as follows. We will describe the language of representatives $\Rep_A^M(\delta)$ of a given piece $\delta \in \Delta$ as the set of ancestors of $\DRep_A^M(\delta)$ under a certain monadic rewriting system, which in turn is controlled by the group of units $U(M)$. Because $\DRep_A^M(\delta)$ can be understood in terms of $U(M)$ by \eqref{Eq:Drep_def}, this gives an understanding of $\Rep_A^M(\delta)$ in terms of $U(M)$. 

We will define the rewriting system 
\begin{equation}\label{Eq:Rewriting_system_RD}
\mathcal{R}_\Delta = \bigcup_{\substack{p \in \Delta \cup \{ \varepsilon \} \\ |p|\leq 1}} \big\{ (W_p \to p) \mid W_p \in \DRep_A^M(p) \big\}.
\end{equation}
In general, this is an infinite rewriting system. Furthermore, it is not in general a complete rewriting system. However, it has the following desirable property: let $\cc$ be a class of languages closed under inverse homomorphism such that $U(M)$ has word problem in $\cc$. Then for every $p \in \Delta \cup \{ \varepsilon \}$ with $|p| \leq 1$ (i.e. for every right-hand side in $\mathcal{R}_\Delta$), the language $\DRep_A^M(p)$ is in $\cc$, as $\Rep_{X}^{U(M)}(\phi(p)) \in \cc$. Thus $\mathcal{R}_\Delta$ is a monadic $\cc$-rewriting system. For example, in the bicyclic monoid $B$ as above, the set of left-hand sides in $\mathcal{R}_\Delta$ of $p \equiv \varepsilon$ is the language $(bc)^\ast$, which is a regular language, as $U(B)$ has regular word problem. 

Before showing the key technical lemma (Lemma~\ref{Lem:M_admits_pres_s.t._REP=<DeltaRep>}), we introduce some useful terminology. For any word in $\dD^\ast$, we can obtain a word in $\Delta^\ast$ by successively removing left-hand sides of rules in $S(M)$, replacing them by their corresponding right-hand sides. We will consider this process in reverse, attributing the terminology of this idea to Cain \& Maltcev \cite{Cain2014}. First, let $u \in \Delta^\ast$ and factorise $u \equiv \delta_1 \delta_2 \cdots \delta_n$ uniquely, where $\delta_i \in \Delta$ for $1 \leq i \leq n$. Then every non-empty subword of $u$ of the form $\delta_j \delta_{j+1} \cdots \delta_{\ell}$ for $1 \leq j \leq \ell \leq n$ is called a \textit{depth-$0$ inserted word}. Inductively, for $\mu \geq 0$, we define a depth-$(\mu+1)$ insertion as follows: if (1) a right-hand side $s_2$ of a rule $(s_1 \to s_2) \in S(M)$ appearing as a proper non-prefix, non-suffix subword of some depth-$\mu$ inserted word $D$, with $D \in \Delta^\ast$, is replaced by $s_1$, then we call that occurrence of $s_1$ a \textit{depth-$(\mu+1)$ inserted word}, and the reversed rewriting $(s_2 \to s_1)$ is then called a depth-$(\mu+1)$ \textit{insertion}; but (2) if instead the specified occurrence of $s_2$ is a depth-$\mu$ inserted word $D \in \Delta^\ast$, or if $s_2 \equiv \varepsilon$ and does not satisfy the condition in (1), then the word $s_1$ is a depth-$\mu$ inserted word. The reversed rewriting $(s_2 \to s_1)$ is then called a depth-$\mu$ insertion.

We give a concrete example. If $\Delta = \{ d, b, abc\}$, and (for simplicity) we have the rewriting system $\mathcal{T}$ with the rules $\{dbd \to b, abc \to \varepsilon \}$, then an ancestor of the word $u \equiv (abc)(abc)(b) \in \Delta^\ast$ modulo $\mathcal{T}$ might look like:
\begin{equation*}\label{Eq:Insertions}
  u' \equiv (adbdc)(abababccc)(dababccdbdd) \equiv (\underbrace{adbdc}_{\textnormal{depth $0$}}) (ab
  \underbrace{ab \overbrace{abc}^{\textnormal{depth $2$}}c }_{\textnormal{depth $1$}} c)(\underbrace{dab\overbrace{abc}^{\textnormal{depth $1$}}cdbdd}_{\textnormal{depth $0$}}).
\end{equation*}
Thus, the word $abc$ in the middle of the word $u'$ is a depth-$2$ inserted word, and the rewriting of the leftmost term $abdbc$ to $abc$ is via the reverse of a depth-$0$ insertion $(b \to dbd)$. Now, just as in \cite[Example~4.2]{Cain2014}, it is clear by definition of insertions (using no properties of the rewriting systems involved) that since $u' \xra{\mathcal{T}} u$, we can perform this rewriting by first rewriting the depth-$2$ insertions in reverse, then the depth-$1$ insertions in reverse, and finally have a depth-$0$ inserted word in $\Delta^\ast$, which is then rewritten to $u$.

\begin{lemma}\label{Lem:M_admits_pres_s.t._REP=<DeltaRep>}
Let $M$ be a finitely presented special monoid. Then $M$ admits a special presentation, with finite generating set $A$ and minimal invertible pieces $\Delta$, such that for every $\delta \in \Delta$ we have $\Rep_A^M(\delta) = \langle \DRep_A^M(\delta) \rangle_{\mathcal{R}_\Delta}$. 
\end{lemma}
\begin{proof}
Let $M$ be as given. By Proposition~\ref{Prop:Every_admits_small_subpiece}, $M$ admits a special presentation satisfying the small subpiece condition. Thus, let us assume that $M$ is given by such a presentation
\begin{equation}\label{Eq:A_pres_satisfying_small_subpiece}
\pres{Mon}{A}{w_1 = 1, w_2 = 1, \dots, w_k=1}
\end{equation}
with pieces $\Delta$. We will prove that for this presentation $\Rep_A^M(\delta) = \langle \DRep_A^M(\delta) \rangle_{\mathcal{R}_\Delta}$.

$(\supseteq)$ Let $w \in \DRep_A^M(\delta)$ and $w' \in A^\ast$ be arbitrary words such that $w' \in \langle w \rangle_{\mathcal{R}_\Delta}$, i.e. $w' \xra{\mathcal{R}_\Delta} w$. Now, for every rule $(W_p, p) \in \mathcal{R}_\Delta$, we have by definition that $p =_M W_p$. Thus, by induction on the number of rules applied in rewriting $w'$ to $w$, we have $w' =_M w$. As $w =_M \delta$, we have $w' \in \Rep_A^M(\delta)$.  

$(\subseteq)$ Let $w \in \Rep_A^M(\delta)$. Then $w =_M \delta$, so $w$ is invertible. In particular, $w \in \dD^\ast$ by Lemma~\ref{Lem:dD_is_invertibles}, and there is some $u \in \Delta^\ast$ such that $w \xra{S} u$. By the earlier reasoning, we can thus obtain $w$ from $u$ by first performing all depth-$0$ insertions, then all depth-$1$ insertions, etc. until after performing a finite number insertions we obtain $w$. Let $\mu \geq 0$ be the highest depth of any such insertion performed. 

We claim that $w \in \langle \DRep_A^M(u)\rangle_{\mathcal{R}_\Delta}$ by induction on this $\mu$. The base case $\mu=0$ is clear, for then $w \in \Delta^\ast$. As for every rule $(s_1 \to s_2) \in S(M)$ we have $s_1 =_M s_2$, it follows by induction on the number of rules applied in rewriting $w \xra{S} u$ that $w =_M u$. As $w \in \Delta^\ast$, we have $w \in \DRep_A^M(u) \subseteq \langle \DRep_A^M(u)\rangle_{\mathcal{R}_\Delta}$. Assume, then, for induction that the claim is true for some $\mu \geq 0$, and that $w$ requires depth-$(\mu+1)$ insertions (but no higher). In the fixed rewriting $w \xra{S} u$, let $u' \in A^\ast$ be such that $w \xr{S} u' \xra{S} u$. Then the rewriting $w \xr{S} u'$ is by replacing the depth-$(\mu+1)$ inserted word $s_1 \in \Delta^\ast$ in $w$ with the word $s_2$, where $s_2$ is a proper non-prefix non-suffix subword of either (I) a depth-$(\mu+1)$ word $Q \in \Delta^\ast$, or (II) a depth-$\mu$ inserted piece $Q \in \Delta$; and where $(s_1 \to s_2) \in S(M)$ is the specified rule. 

In case (I), write $Q \equiv Q_0 s_2 Q_1$, where necessarily $Q_0, Q_1 \in \Delta^\ast$. As $Q$ is a depth-$(\mu+1)$ inserted word in $u'$, the word $u'$ can be obtained from some word $u'' \in A^\ast$ by replacing a depth-$(\mu+1)$ or a depth-$\mu$ inserted word $s_3 \in \Delta^\ast$ in $u''$ with $Q$. That is, there is some rule $(Q \to s_3) \in S(M)$, which rewrites $u' \xr{S} u''$. But as $Q_0 s_1 Q_1 =_M Q_0 s_2 Q_1 \equiv Q =_M s_3$, and $|Q_0 s_1 Q_2| \geq |Q_0 s_2 Q_1| = |Q| \geq |s_3|$, we have $(Q_0 s_1 Q_1 \to s_3) \in S(M)$. Hence, we can obtain $u$ already from $u''$ by performing the insertion of replacing $s_3$ by $s_1$, i.e. $u \xr{S} u''$ by the rule $(Q_0 s_1 Q_1 \to s_3)$, thus reducing the rewriting $w \xra{S} u$ by one step; we may thus by another induction assume without loss of generality that $w$ is obtained from $u'$ as in case (II). 

Thus, assume case (II), i.e. $Q \in \Delta$ is a depth-$\mu$ inserted word in $u'$, and $s_2$ appears as a proper non-suffix non-prefix subword of the piece $Q$. As the presentation satisfies the small subpiece condition, it follows from $s_2 \in \Delta^\ast$ that the subpiece $s_2$ satisfies $|s_2| \leq 1$. Hence also $s_2 \in \Delta \cup \{ \varepsilon \}$. As $s_1 =_M s_2$ and $s_1 \in \Delta^\ast$, we have $s_1 \in \DRep_A^M(s_2)$. Hence, $(s_1 \to s_2) \in \mathcal{R}_\Delta$. Thus, $w$ can be rewritten to $u'$ in a single application of a rule from $\mathcal{R}_\Delta$, and so, by repeating the same step for all depth-$(\mu+1)$ insertions, we find that there is a word $w' \in A^\ast$ such that (1) $w \xra{\mathcal{R}_\Delta} w'$; and (2) $w'$ can be obtained from the word $u$ with at most depth-$\mu$ insertions. By the inductive hypothesis, thus $w' \in \langle u \rangle_{\mathcal{R}_\Delta}$, and hence also $w \in \langle u \rangle_{\mathcal{R}_\Delta}$. Now $u \in \DRep_A^M(\delta)$, as $u \in \Delta^\ast$ and $u =_M \delta$, and so we conclude that $w \in  \langle \DRep_A^M(\delta) \rangle_{\mathcal{R}_\Delta}$, which is what was to be shown. 
\end{proof}

We conclude:

\begin{theorem}
Let $M = \pres{Mon}{A}{w_1 =1, w_2 = 1, \dots, w_k = 1}$. Let the group of units $U(M)$ of $M$ be generated by a finite set $X$, and let $\cc$ be a super-$\AFL$. Then $M$ admits a special presentation, with invertible pieces $\Delta$, such that 
\[
\WP_X^{U(M)} \in \cc \quad \implies \quad \Rep_A^M(\delta) \in \cc \:\: \textnormal{for all $\delta \in \Delta$.}
\]
\end{theorem}
\begin{proof}
By Lemma~\ref{Lem:M_admits_pres_s.t._REP=<DeltaRep>}, $M$ admits a presentation satisfying the conclusions of that lemma. Suppose $M$ is given by such a presentation, with pieces $\Delta$, and $X$ as usual. Suppose $\WP_X^{U(M)} \in \cc$, and let $\delta \in \Delta$ be any piece. Then for every $p \in \Delta \cup \{ \varepsilon \}$ with $|p|\leq 1$, we have $\DRep_A^M(p) \in \cc$. Hence the left-hand side of every letter or the empty word in $\mathcal{R}_\Delta$ is in $\cc$, so $\mathcal{R}_\Delta$ is a monadic $\cc$-rewriting system. The conclusions of Lemma~\ref{Lem:M_admits_pres_s.t._REP=<DeltaRep>} being satisfied yields that $\Rep_A^M(\delta) = \langle \DRep_A^M(\delta) \rangle_{\mathcal{R}_\Delta}$. As $\cc$ has the monadic ancestor property and $\DRep_A^M(\delta) \in \cc$, we find that $\Rep_A^M(\delta) \in \cc$, as required. 
\end{proof}

This theorem is thus a full description of what the invertible words of a special monoid look like (we shall soon see that the assumption on the presentation becomes somewhat unimportant). We shall now use this description to language-theoretically describe when two invertible words are equal in a special monoid. 

\subsection{The invertible word problem} By the normal form lemma Lemma~\ref{Lem:Zhangs_lemma}, understanding equality of invertible words is tantamount to understanding equality of words, i.e. in describing in the word problem. We note that, on the surface, equality of invertible words (over $A^\ast$) is quite distinct from equality of elements in $U(M)$. Motivated by the word problem $\WP_A^M$ for $M$, we define the \textit{invertible word problem} of $M$ (with respect to $A$)  as 
\begin{equation}\label{Eq:InvP_def}
\InvP_A^M = \{ w_1 \# w_2^\trev \mid w_1, w_2 \in A^\ast \textnormal{ invertible, and } w_1 =_M w_2 \}.
\end{equation}
By Lemma~\ref{Lem:dD_is_invertibles}, $\InvP_A^M = \WP_A^M \cap \dD^\ast \# (\dD^\trev)^\ast = \WP_A^M \cap \mathfrak{M}^\ast \# (\mathfrak{M}^\trev)^\ast$. 

We first provide the language-theoretic version of Lemma~\ref{Lem:Zhangs_lemma}. 

\begin{lemma}\label{Lem:InvP=>WP}
Let $\cc$ be a super-$\AFL$. Then $\InvP_A^M \in \cc \implies \WP_A^M \in \cc$.
\end{lemma}
\begin{proof}
The result is obvious using Lemma~\ref{Lem:Zhangs_lemma} and the alternating products introduced by the author \cite{NybergBrodda2020c, NybergBrodda2021f}. We provide a direct proof instead. Let $\mathcal{T}$ be the rewriting system over the alphabet $A \cup \{ \# \}$ with the rules
\begin{equation}\label{Eq:RewritingT(Inv=>WP)}
\{ W \to \# \mid W \in \InvP_A^M \} \cup \{ a \# a \to \# \mid a \in A \}.
\end{equation}
Then $\mathcal{T}$ is a monadic rewriting system. It is obviously a $\cc$-rewriting system, as the language of left-hand sides of the symbol $\#$ in \eqref{Eq:RewritingT(Inv=>WP)} is the union of two languages in $\cc$, and $\cc$ is closed under unions. We claim that $\langle \# \rangle_{\mathcal{T}} = \WP_A^M$. Upon proving this, we will conclude, as $\mathcal{T}$ is a monadic $\cc$-rewriting system, $\cc$ has the monadic ancestor property, and $\{ \# \} \in \cc$, that $\WP_A^M \in \cc$, as desired.

$(\subseteq)$ Let $w \in \langle \# \rangle_{\mathcal{T}}$. Then there exists some least $\mu \geq 0$ and words $u_i \in A^\ast$, where $0 \leq i \leq \mu$ such that 
\begin{equation}\label{Eq:RewriteW->hash}
w \equiv u_0 \xr{\mathcal{T}} u_1 \xr{\mathcal{T}} \cdots \xr{\mathcal{T}} u_{\mu-1} \xr{\mathcal{T}} u_\mu \equiv \#.
\end{equation}
We prove by induction on $\mu$ that $w \in \WP_A^M$. The base case $\mu=0$ implies $w \equiv \#$, and of course $\# \in \WP_A^M$. Suppose $\mu > 0$ and that the claim is true for all rewritings of the form \eqref{Eq:RewriteW->hash} requiring at most $\mu-1$ steps. Now $w_1 \xr{\mathcal{T}}^{\mu-1} \#$, so by the inductive hypothesis we have $w_1 \in \WP_A^M$. Thus, we can write $w_1 \equiv w_1' \# (w_1'')^\trev$, where $w_1', w_1'' \in A^\ast$ satisfy $w_1' =_M w_1''$. Let $(r,s) \in \mathcal{T}$ be the rule which rewrites $w_0 \xr{\mathcal{T}} w_1$. Then $s \equiv \#$ by \eqref{Eq:RewritingT(Inv=>WP)}, so $w \equiv w_0 \equiv w_1' r (w_1'')^\trev$. As $r$ is a left-hand side of a rule in $\mathcal{T}$, thus $s \equiv u \# v^\trev$, where either $u \# v^\trev \in \InvP_A^M$, or else $u \equiv a \equiv v$, where $a \in A$. In either case, $u =_M v$. Thus 
\[
w \equiv w_1' r (w_1'')^\trev \equiv w_1' u \# v^\trev (w_1'')^\trev \equiv (w_1'u) \# (w_1'' v)^\trev,
\]
and as $w_1' u =_M w_1'' u =_M w_1'' v$, thus $w \in \WP_A^M$. 

$(\supseteq)$ Suppose that $w \in \WP_A^M$. Then $w \equiv u \# v^\trev$ for some $u, v \in A^\ast$ with $u=_M v$. By Lemma~\ref{Lem:Zhangs_lemma}, we can factorise $u$ and $v$ uniquely as 
\begin{align*}
u \equiv u_0a_1 u_1 \cdots a_m u_m, \quad v \equiv v_0a_1 v_1 \cdots a_m v_m, 
\end{align*}
respectively, where for every $0 \leq i \leq m$ we have $a_i \in A$, $u_i =_M v_i$, and $u_i$ (resp. $v_i$) is a maximal invertible factor of $u$ (resp. $v$). We prove the claim by induction on this $m$. The base case $m=0$ is clear, for then $u \equiv u_0, v \equiv v_0$ are invertible, and thus $w \equiv u \# v^\trev \in \InvP_A^M$. Hence $(w \to \#) \in \mathcal{T}$, so $w \in \langle \# \rangle_{\mathcal{T}}$. Assume $m>0$ and that the claim holds for $m-1$. As $u_m, v_m$ are invertible, we have $(u_m \# v_m^\trev \to \#), (a_m \# a_m \to \#) \in \mathcal{T}$. Thus 
\begin{align*}
w \equiv u \# v^\trev \equiv u_0a_1 u_1 \cdots u_{m-1} a_m u_m &\# v_m^\trev a_m v_{m-1}^\trev \cdots a_1 v_0^\trev \\ \xra{\mathcal{T}} u_0a_1 u_1 \cdots u_{m-1} &\# v_{m-1}^\trev \cdots a_1 v_0^\trev.
\end{align*}
By the inductive hypothesis, $u_0a_1 u_1 \cdots u_{m-1} \# v_{m-1}^\trev \cdots a_1 v_0^\trev \in \langle \# \rangle_{\mathcal{T}}$, and so also $w \in \langle \# \rangle_{\mathcal{T}}$, as was to be shown.
\end{proof}

Lemma~\ref{Lem:InvP=>WP} shows that understanding $\InvP_A^M$ translates to understanding $\WP_A^M$. We now describe $\InvP_A^M$ in terms of the word problem for $U(M)$.  

\begin{lemma}\label{Lem:UM CF iff INVP CF}
Let $\cc$ be a super-$\AFL$ closed under reversal. Then the special monoid $M$ admits a special monoid presentation, with generators $A$, such that $U(M)$ has word problem in $\cc$ if and only if $\InvP_A^M \in \cc$. 
\end{lemma}
\begin{proof}
By Lemma~\ref{Lem:M_admits_pres_s.t._REP=<DeltaRep>}, $M$ admits a special monoid presentation satisfying the conclusions of that lemma; let $M$ be given by such a presentation, and let $A, \Delta, X$ and $\phi \colon \Delta^\ast \to X^\ast$ be as usual. 

$(\impliedby)$ For notational brevity, write $\Delta_r = \Delta^\trev$. The language $\Delta^\ast \# \Delta_r^\ast$ is a regular language. Let $K = \InvP_A^M \cap (\Delta^\ast \# \Delta_r^\ast)$. As $\cc$ is closed under intersection with regular languages, $K \in \cc$. Now $K$ consists of precisely the words of the form $u \# v^\trev$ where $u, v \in \Delta^\ast$ and $u =_M v$. That is, $K = \WP_\Delta^{U(M)}$, cf. \eqref{Eq: WPDeltaUM}. Hence $\WP_\Delta^{U(M)} \in \cc$, so $U(M)$ has word problem in $\cc$. 

$(\implies)$ For every $\delta \in \Delta$, let $\hs_\delta, \ths_\delta$ be new symbols. Define $\hs_\Delta = \{ \hs_\delta \mid \delta \in \Delta \}$, and $\ths_\Delta = \{ \ths_\delta \mid \delta \in \Delta \}$ such that $\hs_\Delta \cap \ths_\Delta = \varnothing$. Let $R_\delta, R_\delta^r$ be the rewriting systems on the alphabet $A \cup \hs_\Delta \cup \ths_\Delta$ defined by 
\begin{align*}
\cR_\delta := &\bigcup_{\delta \in \Delta} \big\{ (w, \hs_\delta) \mid w \in \Rep_A^M(\delta) \big\}, \\
\cR_\delta^r := &\bigcup_{\delta \in \Delta} \big\{ (w^\trev, \ths_\delta) \mid w \in \Rep_A^M(\delta) \big\}.
\end{align*}
The conclusions of Lemma~\ref{Lem:M_admits_pres_s.t._REP=<DeltaRep>} being satisfied, $\Rep_A^M(\delta) \in \cc$ for every $\delta \in \Delta$; as $\Delta$ is finite and $\cc$ is closed under finite unions, $\cR_\delta$ is a monadic $\cc$-rewriting system. As $\cc$ is closed under reversal, $\cR_\delta^r$ is also a monadic $\cc$-rewriting system, and so the system $\cR_0 := \cR_\delta \cup \cR_\delta^r$ is, too. 

There exists a surjective homomorphism $\varrho \colon (\hs_\Delta \cup \ths_\Delta)^\ast \to \Delta^\ast$ defined, for $\delta \in \Delta$ by $\varrho(\hs_\delta) = \varrho(\ths_\delta) = \delta$. Thus $\hs_\Delta \cup \ths_\Delta$ is a finite generating set for $U(M)$. Let 
\begin{equation}\label{Eq:WPhsths}
L := \WP_{\hs_\Delta \cup \ths_\Delta}^{U(M)} \cap (\hs_\Delta^\ast \# \ths_\Delta^\ast).
\end{equation}
As $U(M)$ has word problem in $\cc$, and $\cc$ is a super-$\AFL$, we have $L \in \cc$. We claim that $\InvP_A^M = \langle L \rangle_{\cR_0} \cap A^\ast \# A^\ast$. As $\cc$ has the monadic ancestor property, this would imply $\langle L \rangle_{\cR_0} \in \cc$ and consequently $\InvP_A^M \in \cc$, completing the proof. 

$(\subseteq)$ Let $w \in \InvP_A^M$ be arbitrary. By Lemma~\ref{Lem:dD_is_invertibles}, there exist $w_1, w_2 \in \dD^\ast$ such that $w \equiv w_1 \# w_2^\trev$ with $w_1 =_M w_2$. Write $w_1 \equiv \vartheta_0 \vartheta_1 \cdots \vartheta_n$ and $w_2 \equiv \vartheta_0' \vartheta_1' \cdots \vartheta'_m$, where $\vartheta_i, \vartheta_j' \in \dD$ for $0 \leq i \leq n$ and $0 \leq j \leq m$. By definition of $\dD$, for every $\vartheta_i, \vartheta_j'$ there exist $\delta_i, \delta_j' \in \Delta$ such that $\vartheta_i =_M \delta_i$ and $\vartheta_j' =_M \delta_j'$, i.e. $\vartheta_i \in \Rep_A^M(\delta_i)$ and $\vartheta_j' \in \Rep_A^M(\delta_j')$. Thus for every $i, j$ as above, $(\vartheta_{\delta_i}, \hs_{\delta_i}), ( (\vartheta_j')^\trev, \ths_{\delta_j'}) \in \cR_0$. Let $W$ be the word $\hs_{\delta_0} \hs_{\delta_1} \cdots \hs_{\delta_n} \# \ths_{\delta'_m} \ths_{\delta'_{m-1}} \cdots \ths_{\delta'_0}$. Then $w \xra{\mathcal{R}_0} W$. As $\delta_0 \delta_1 \cdots \delta_n =_M \delta_0' \delta_1' \cdots \delta_m'$, it follows that $W \in \WP_{\hs_\Delta \cup \ths_\Delta}^{U(M)}$, and as $W \in \hs_\Delta^\ast \# \ths_\Delta^\ast$, we have $W \in L$. We conclude that $w \in \langle L \rangle_{\cR_0}$. As $w \in A^\ast \# A^\ast$, thus $\langle L \rangle_{\cR_0} \cap A^\ast \# A^\ast$. 

$(\supseteq)$ Let $w \in \langle L \rangle_{\cR_0} \cap A^\ast \# A^\ast$. Let $W \in L$ be such that $w \xr{\cR_0} W$. Then there exist $\delta_i, \delta_j' \in \Delta$ and $\hs_{\delta_i} \in \hs_\Delta, \ths_{\delta_j}' \in \ths_\Delta$ for $0 \leq i \leq n$ and $0 \leq j \leq m$ such that 
\begin{equation}
W \equiv \hs_{\delta_0}  \hs_{\delta_1} \cdots \hs_{\delta_n} \# \ths_{\delta_m'} \ths_{\delta_{m-1}'} \cdots \ths_{\delta_0'},
\end{equation}
with $\delta_0 \delta_1 \cdots \delta_n =_M \delta_0' \delta_1' \cdots \delta_m'$. Let $u \in (A \cup \hs_\Delta \cup \ths_\Delta \cup \{ \# \})^\ast$ be such that $u \xra{\cR} W$. As no left-hand side of a rule in either $\cR_\delta$ or $\cR_\delta^r$ contains an occurrence of any letter from $\hs_\Delta$ or $\ths_\Delta$, it easily follows that $u$ can be written as $\alpha_0 \alpha_1 \cdots \alpha_n \# \alpha_m' \alpha_{m-1}' \cdots \alpha_0'$, where for $0 \leq i \leq n$, $\alpha_i$ is either (i) $\hs_{\delta_i}$, or else (ii) $u_i$, where $u_{i} \in \Rep_A^M(\delta_i)$; and similarly for $0 \leq j \leq m$, $\alpha_j'$ is either (i') $\ths_{\delta_j'}$, or else (ii') $v_j^\trev$, where $v_j \in \Rep_A^M(\delta_j')$. In particular, as $w \in A^\ast \# A^\ast$ and $w \xr{\cR} W$, there exist $u_i, v_j$ as above such that 
\[
w \equiv u_0 u_1 \cdots u_n \# v_m^\trev v_{m-1}^\trev \cdots v_0^\trev \equiv  u_0 u_1 \cdots u_n \# (v_0 v_1 \cdots v_m)^\trev,
\]
with $u_i \in \Rep_A^M(\delta_i)$ and $v_j \in \Rep_A^M(\delta_j')$ for all $0 \leq i \leq n$ and $0 \leq j \leq m$. As $\delta_0 \delta_1 \cdots \delta_n =_M \delta_0' \delta_1' \cdots \delta_m'$, it it follows that $u_0 u_1 \cdots u_n =_M v_0 v_1 \cdots v_m$. Furthermore, every $u_i$ and every $v_j$ is invertible. Thus $w \in \InvP_A^M$. 
\end{proof}

\begin{example}
We give an example of the system $\cR_0$ from the proof of Lemma~\ref{Lem:UM CF iff INVP CF} for a concrete special monoid. Let $M_4 = \pres{Mon}{a,b,c}{abc=1, b^2 =1}$. Then it is easily seen that $\Delta = \{ abc, b\}$. Thus $\hs_\Delta = \{ \hs_{abc}, \hs_b \}$ and $\ths_\Delta = \{ \ths_{abc}, \ths_b \}$. Note that $ab^{17}c =_{M_4} abc$, though certainly $ab^{17}c \not\in \Delta$. Thus we have the rule $(ab^{17}c \to \hs_{abc})$ in $\cR_\delta \subset \cR_0$. Similarly, $bab^3 cb^2 =_{M_4} b$, so $(b^2cb^3ab \to \ths_{b}) \in \cR_\delta^r \subset \cR_0$. 
\end{example}

Let $\cc$ be closed under inverse homomorphism. While the word problem for $M$ being in $\cc$ is invariant under change of generating set, it is not obvious that the \textit{invertible} word problem for $M$ being in $\cc$ is also such an invariant. We suspect that this is the case, but have been unable to find a proof. In any case, Lemma~\ref{Lem:UM CF iff INVP CF} provides a sufficiently good description of $\InvP_A^M$ in terms of $U(M)$ for our purposes.

\section{The word problem}

\noindent We have now completely described the invertible word problem for $M$ in terms of $U(M)$. This yields a proof of the main theorem of this article. 

\begin{theoremB}\label{Thm:Main_thm}
Let $M$ be a finitely presented special monoid. Let $\cc$ be a super-$\AFL$ closed under reversal. Then $M$ has word problem in $\cc$ if and only if the group of units $U(M)$ of $M$ has word problem in $\cc$. 
\end{theoremB}
\begin{proof}
$(\implies)$ As $M$ is finitely presented, $U(M)$ is finitely generated, and if $M$ has word problem in $\cc$, then by \cite[Proposition~8(a)]{Hoffmann2012} as $\cc$ is a super-$\AFL$ any finitely generated submonoid of $M$ has word problem in $\cc$. Thus $U(M)$ has word problem in $\cc$. 

$(\impliedby)$ By Lemma~\ref{Lem:UM CF iff INVP CF}, $M$ admits a presentation satisfying the conclusions of Lemma~\ref{Lem:UM CF iff INVP CF}. Let $A$ be the finite generating set for $M$ in this presentation, and let further $\Delta, X$ be the minimal pieces and associated set. Then, as $\cc$ is closed under inverse homomorphism, $\WP_{X'}^{U(M)} \in \cc$. As the conclusions of Lemma~\ref{Lem:UM CF iff INVP CF} are satisfied, we have $\InvP_A^M \in \cc$. By Lemma~\ref{Lem:InvP=>WP}, thus $\WP_A^M \in \cc$. As $\cc$ is closed under inverse homomorphism, the word problem for $M$ (with respect to any finite generating set) is thus also in $\cc$. 
\end{proof}

We remark that the assumption of ``finitely presented'' cannot in general be dropped from the statement of Theorem~\ref{Thm:Main_thm}, see Remark~\ref{Rem:Bad_special_monoids}(3). As one application of Theorem~\ref{Thm:Main_thm}, we can take the class $\IND$ of \textit{indexed} languages, which is a super-$\AFL$ closed under reversal (see \cite{Engelfriet1985}). Furthermore, as the class $\CF$ of context-free languages is a super-$\AFL$ closed under reversal, we find the following.

\begin{corollary}\label{Cor:MS}
A finitely presented special monoid has context-free word problem if and only if its group of units is virtually free. 
\end{corollary}

As before, the assumption of ``finitely presented'' in Corollary~\ref{Cor:MS} cannot in general be dropped, see Remark~\ref{Rem:Bad_special_monoids}(3). Now, as every group is a special monoid (equal to its own group of units), this is a full generalisation of the Muller-Schupp theorem. In 2004, Duncan \& Gilman asked for a characterisation of monoids with context-free word problem \cite[Question~4]{Duncan2004}. We have thus answered this question completely for the class of special monoids. Hoffmann et al. \cite[p. 97]{Hoffmann2012} write ``the depth of the Muller-Schupp result and its reliance on the geometrical structure of Cayley graphs of groups suggests that a generalization to semigroups could be very hard to obtain''. Nevertheless, the above Corollary~\ref{Cor:MS} is a generalisation of this sort, free from any geometric structure. We note that every monoid with context-free word problem is \textit{word-hyperbolic} (in the sense of \cite{Duncan2004}, see also \cite{Cain2012}); thus a further corollary is that a finitely presented special monoid with virtually free group of units is word-hyperbolic. For example, $\pres{Mon}{a,b,c,d}{abcdab=1}$ is word-hyperbolic, as it has infinite cyclic group of units. 

\begin{remark}\label{Rem:Bad_special_monoids}
We give some examples illustrating how context-free special monoids, despite the similarities suggested by Theorem~\ref{Thm:Main_thm}, differ from groups. All the following statements are false if one substitutes ``group'' for ``special monoid''. 
\begin{enumerate}[label=(\arabic*)]
\item There exists a (finitely generated) context-free special monoid which cannot be finitely presented. Namely, let $M_5 = \pres{Mon}{a,b,c}{ab^ic = 1 \: (i \geq 1)}$. Then $M_5$, being defined by a context-free complete monadic rewriting system, has context-free word problem \cite[Corollary~3.8]{Book1982b}, but obviously $M_5$ cannot be finitely presented. By contrast, every context-free group is finitely presented.
\item There exists a finitely generated special monoid with context-free word problem, whose group of units is not finitely generated (indeed none of whose maximal subgroups is finitely generated). We refer the reader to Nyberg-Brodda \cite{NybergBrodda2021d} for this example, which answered in the negative a question of Brough, Cain \& Pfeiffer.
\item There exists a finitely generated special monoid with context-free group of units, but such that the word problem for the monoid is not context-free. Namely, let $M_6 = \pres{Mon}{a,b,c,d}{ab^ic^jd = 1 \: (i, j \geq 2, j = i^2)}$. Write, for brevity,  $A = \{ a, b,c d\}$. Then the set of defining relations of $M_6$ clearly forms a complete monadic rewriting system defining $M_6$. For any defining word $w_{i,j}: \equiv ab^ic^jd$, as no rule of the complete system begins with $b, c$ or $d$, no word $bu, cu$, or $du$ is (right) invertible for any $u \in A^\ast$, and hence the factorisation of $w_{i,j}$ into minimal invertible factors is trivial. As the minimal invertible factors generates $U(M)$ (see e.g. Gray \& Ru\v{s}kuc \cite{Gray2021}), we conclude that $U(M_6) = 1$. On the other hand, using the complete rewriting system it is straightforward to see that
\[
\Rep_{A}^{M_6}(1) \cap ab^\ast c^\ast d = \{ ab^i c^j d \mid i, j \geq 2, j = i^2 \},
\]
but this right-hand side is well-known to not be a context-free language. Thus $\Rep_A^{M_6}(1)$ cannot possibly be context-free, and so $M_6$ cannot have context-free word problem, see Theorem~\ref{Thm:Equiv_of_three_notions}(i)$\implies$(iii) (this direction does not require finite presentability). 
\end{enumerate}
\end{remark}

Despite the examples given in Remark~\ref{Rem:Bad_special_monoids}, groups and special monoids share many language-theoretic similarities. We present these now, all of which essentially follow from Theorem~\ref{Thm:Main_thm}. 

\subsection{Representatives of words}

Recall that to speak of a ``context-free group'' or a ``regular group'' means to speak of the language words representing the identity element of this group. For monoids, such a language is not, in general, interesting; if all defining relations of a monoid $\Pi$ are of the form $u = v$ with $u, v \in A^+$ non-empty, for example, then only the empty word represents the identity element. On the other hand, for special monoids, this language completely characterises the language-theoretic behaviour of the monoid:

\begin{theorem}\label{Thm:Equiv_of_three_notions}
Let $M = \pres{Mon}{A}{w_1 = 1, w_2 = 1, \dots, w_k = 1}$ be a finitely presented special monoid. Let $\cc$ be a super-$\AFL$ closed under reversal. Then the following are equivalent:
\begin{enumerate}[label=(\roman*), font=\normalfont]
\item $M$ has word problem in $\cc$.
\item For every word $w \in A^\ast$, $\Rep_A^M(w) \in \cc$.
\item $\IP_A^M = \{ w \mid w =_M 1, w \in A^\ast\}$ is in $\cc$.
\end{enumerate}
\end{theorem}
\begin{proof}
(i)$\implies$(ii). Let $w \in A^\ast$. If $\WP_A^M \in \cc$, then the right quotient $\WP_A^M  / \{ \# w^\trev \}$ is in $\cc$, as $\cc$ is an $\AFL$ and hence closed under right quotients with regular languages. But this quotient language is just $\{ u \in A^\ast \mid u =_M w \} = \Rep_A^M(w)$. 

(ii)$\implies$(iii) Obvious.

(iii)$\implies$(i) Suppose $\IP_A^M \in \cc$. Let $\Delta$ be the pieces of the presentation. Then, as $\cc$ is closed under intersection with regular languages, $\IP_A^M \cap \Delta^\ast \in \cc$. But 
\[
\IP_A^M \cap \Delta^\ast = \{ w \in \Delta^\ast \mid w =_M 1\} =  \IP_{\Delta}^{U(M)}
\]
where $\Delta$ is considered as a generating set for $U(M)$ (cf. \eqref{Eq: WPDeltaUM}). As $U(M)$ is a group, having $\IP_\Delta^{U(M)} \in \cc$ implies by \cite[Theorem~5.3]{Duncan2004} that $U(M)$ has word problem in $\cc$. Thus, by Theorem~\ref{Thm:Main_thm}, $M$ has word problem in $\cc$.
\end{proof}

Taking $\cc = \CF$, the class of context-free languages, Theorem~\ref{Thm:Equiv_of_three_notions} answers a 1992 question first posed by Zhang \cite[Problem~1]{Zhang1992}, who asked: if the group of units $U(M)$ of a special monoid is context-free, are the languages $\Rep_A^M(w)$ for $w \in A^\ast$ context-free? By Theorem~\ref{Thm:Main_thm} and (i)$\implies$(ii) in Theorem~\ref{Thm:Equiv_of_three_notions}, the answer to Zhang's question is thus affirmative. For $w \in A^\ast$, the set $\Rep_A^M(w)$ is called a \textit{basic congruential language} by Zhang. The study of congruences on the free monoid all of whose congruence classes are context-free languages has a long history, which we do not expand on here; we mention only the early work by Cochet \cite{Cochet1971}.
 
\begin{corollary}
Let $M, \cc$ be as in Theorem~\ref{Thm:Equiv_of_three_notions}. If the group of units $U(M)$ has word problem in $\cc$, then $\mathfrak{M}^\ast = \dD^\ast \in \cc$. 
\end{corollary} 
\begin{proof}
Obvious by Theorem~\ref{Thm:Main_thm} and Theorem~\ref{Thm:Equiv_of_three_notions}, as $\dD^\ast$ is the right quotient of $\WP_A^M$ (which is in $\cc$) by the regular language $(\Delta^\trev)^\ast$, and $\cc$ is an $\AFL$.
\end{proof}

In particular, if the group of units of a (finitely presented) special monoid $M$ (generated by $A$) is virtually free, then the set of all invertible words is a context-free language. It is thus decidable in sub-cubic time (see the proof of Proposition~\ref{Prop:Solve_in_n23728639}) whether a given word $w \in A^\ast$ represents an invertible element of $M$ or not. 

Finally, we remark that in the proof of Theorem~\ref{Thm:Equiv_of_three_notions}, it is easy use Lemma~\ref{Lem:Zhangs_lemma} to allow for replacing the empty word $1$ in Theorem~\ref{Thm:Equiv_of_three_notions}(iii) with \textit{any} word $w \in A^\ast$. That is, $\Rep_A^M(w)$ is in $\cc$ for \textit{any} word $w \in A^\ast$ if and only if it is in $\cc$ for \textit{every} word $w \in A^\ast$. We leave the details to the reader. Theorem~\ref{Thm:Equiv_of_three_notions} shows that special monoids and groups share similar language-theoretic behaviour, and that for special monoids the definition of the word problem as given here is a good generalisation of the word problem for groups. 

\subsection{Some decision problems}

In 1992, Zhang \cite[Problem~3]{Zhang1992d} asked: given a special one-relation monoid $\pres{Mon}{A}{w=1}$, is it decidable whether the congruence class of every word is a context-free language? The below theorem, when combined with Theorem~\ref{Thm:Equiv_of_three_notions}, gives a complete and affirmative answer to Zhang's question.

\begin{theoremB}\label{Thm:Dec_if_w=1_CF}
It is decidable whether a special one-relation monoid $\pres{Mon}{A}{w=1}$ has context-free word problem. 
\end{theoremB}
\begin{proof}
Let $M = \pres{Mon}{A}{w=1}$ be a special one-relation monoid. By Theorem~\ref{Thm:Main_thm}, it suffices to decide whether the group of units $U(M)$ is virtually free. Now $U(M)$ is a one-relator group $\pres{Gp}{X}{r=1}$, where $X$ and $r$ are effectively computable from the given presentation of $M$ by Adian's overlap algorithm. First, we decide whether $U(M)$ is torsion-free: this is the case if and only if $r$ is not equal in the free group $F_X$ on $X$ to a proper power $u^n$ of some other word $u$, $n >1$. If $U(M)$ is torsion-free, then $U(M)$ is virtually free if and only if it is free \cite{Stallings1968}, and $U(M)$ is in turn free if and only if $r$ is empty or a primitive element $F_X$ by \cite[Theorem~4]{Whitehead1936}. As primitivity of $r$ can be decided by Whitehead's algorithm \cite{Whitehead1936}, this yields the result if $U(M)$ is torsion-free. If $U(M)$ is not torsion-free, then we can (uniquely) write $r = u^n$ as above, with equality in $F_X$. Then $U(M)$ is virtually free if and only if $u$ is primitive in $F_X$ by \cite[Theorem~3]{Fischer1972}, which can again be decided by Whitehead's algorithm.
\end{proof}

We remark that Zhang additionally asked if this problem can be solved in polynomial time; the above algorithm is polynomial-time, thus also answering this part of the question affirmatively. Indeed, the only part of the algorithm which is not obviously polynomial-time is using the Whitehead algorithm to check if a word is primitive; but this can be done in quadratic time, see e.g. \cite{Myasnikov2003}.

\begin{example}
As an illustration of Theorem~\ref{Thm:Dec_if_w=1_CF}, we give two examples. 

\begin{enumerate}[label=(\arabic*)]
\item Let $M_7 = \pres{Mon}{a,b}{(abaabbab)^2 = 1}$. Then by Adian's overlap algorithm the defining word factors into minimal invertible pieces as $((ab)(aabb)(ab))^2$, and the group of units of $M_7$ is isomorphic to $\pres{Gp}{x_1,x_2}{(x_1x_2x_1)^2 = 1}$. The word $x_1x_2x_1$ is readily seen to be a primitive word; we may use the Nielsen transformations $x_2 \mapsto x_2x_1^{-1}$ followed by $x_1 \mapsto x_1^{-1}x_2^{-1}$ to transform $x_1x_2x_1$ into $x_1 x_2$, followed by $x_1$. Thus $U(M_7) \cong \pres{Gp}{x_1,x_2}{x_1^2 =1}$, so $U(M_7) \cong C_2 \ast \mathbb{Z}$ is virtually free (explicitly, the subgroup generated by $x_2$ and $x_1 x_2 x_1^{-1}$ has index $2$ in $U(M_7)$ and is free of rank $2$). We conclude that $M_7$ has context-free word problem. 
\item Let $M_8 = \pres{Mon}{a,b,c}{acabcabcac=1}$. We factor, as before, the defining word into minimal invertible pieces as $(ac)(abc)(abc)(ac)$, which yields us $U(M_8) \cong \pres{Gp}{x_1,x_2}{x_1 x_2^2 x_1 = 1}$. However, $x_1x_2^2 x_1$ is not a primitive word; indeed, $U(M_8) \cong \pres{Gp}{x_1,x_2}{x_1^2x_2^2 =1}$, the fundamental group of the Klein bottle, which is virtually $\mathbb{Z}^2$, and therefore not even hyperbolic. In particular, $U(M_8)$ is \textit{not} virtually free by Theorem~\ref{Thm:Main_thm}, so $M_8$ is not context-free.
\end{enumerate}
\end{example} 

Theorem~\ref{Thm:Dec_if_w=1_CF} can be seen as partial progress in studying the following question, which does not appear to have been studied anywhere previously.

\begin{question}\label{Quest:Dec_if_OR_CF?}
Let $M = \pres{Mon}{A}{u=v}$ be a one-relation monoid. Is it decidable whether $M$ has context-free word problem?
\end{question}

Note that whereas it is still an open problem whether the word problem for all one-relation monoids is decidable (see Nyberg-Brodda \cite{NybergBrodda2021b} for a survey of this problem), this does not preclude classifying those one-relation monoids with context-free word problem; cf. e.g. the classification by Shneerson \cite{Shneerson1972a, Shneerson1972b} of the one-relation monoids which satisfy some non-trivial identity. For further progress on Question~\ref{Quest:Dec_if_OR_CF?}, see Nyberg-Brodda \cite{NybergBrodda2020c}.

We turn to the word problem for special monoids, considered as a decision problem. Let $M$ be a finitely presented special monoid, generated by $A$ and whose group of units has decidable word problem. The method devised by Makanin \cite{Makanin1966} or indeed Zhang \cite{Zhang1992} to decide whether, for $u, v \in A^\ast$, we have $u =_M v$, is exponential in $f(|u|+|v|)$, where $f$ is the complexity of the word problem for the group of units. It would be interesting to see to what extent this can be sharpened. The best we can do for now is the following, which follows directly from Theorem~\ref{Thm:Main_thm}.

\begin{proposition}\label{Prop:Solve_in_n23728639}
Let $M$ be a finitely presented special monoid, finitely generated by $A$, with virtually free group of units. Then the word problem for $M$, with input $u, v \in A^\ast$, is decidable in $O(n^{2.3728639})$-time, where $n = |u|+|v|$.
\end{proposition}
\begin{proof}
By Theorem~\ref{Thm:Main_thm}, there exists a context-free grammar $\Gamma_M$ generating $\WP_A^M$. Thus the word problem reduces to checking whether $u \# v^\trev \in \cl(\Gamma_M)$. By a result of Valiant \cite{Valiant1975}, checking membership for a word of length $n$ in the language of a context-free grammar reduces to the problem of multiplying $n \times n$-matrices with entries in $\operatorname{GF}(2)$. The current best known algorithm for this latter problem is $O(n^{2.3728639})$, due to Le Gall \cite{LeGall2014}, see also Williams \cite{Williams2012}.
\end{proof}

We mentioned earlier that any context-free monoid is word-hyperbolic. Unlike for hyperbolic groups (for which the word problem is decidable in linear time), the best known result along these lines for word-hyperbolic monoids is that the word problem can be solved in polynomial time \cite{Cain2016}. No upper bound on the degree for the polynomials that arise in this way is known to exist, but the best current known algorithm cannot be faster than $O(n^5 \log n)$. 

\begin{conjecture}\label{Conj:VF=>quadratic}
Let $M$ be a finitely presented special monoid with virtually free group of units. Then the word problem is decidable in quadratic time. 
\end{conjecture}

The conjectured fastest time for matrix multiplication is $O(n^2)$; if that conjecture is true, then the proof of Proposition~\ref{Prop:Solve_in_n23728639} demonstrates that Conjecture~\ref{Conj:VF=>quadratic} is also true. In line with this conjecture, we ask the following rather broad question. 

\begin{question}
Let $M$ be a finitely presented special monoid such that the word problem for the group of units of $M$ is decidable in $O(f(n))$ for some function $f$. Does there always exist a polynomial $g$ such that the word problem for $M$ is decidable in $O(g(f(n))$?
\end{question}

Finally, we turn to the rational subset membership problem. A subset $K \subseteq M$ is \textit{rational} if and only if there exists a regular language $L \subseteq A^\ast$ such that $K = \pi(L)$ (see \cite{Eilenberg1969} for this definition). The rational subset membership problem is said to be decidable if, given as input a word $w \in A^\ast$ and a regular language $K \subseteq A^\ast$ (given e.g. as a finite-state automaton accepting the language), then we can decide if $\pi(w) \in \pi(K)$. Every singleton element of a monoid $M$ is a rational subset of $M$; as is every finitely generated submonoid; as is, for any $u \in A^\ast$, the set of elements $m \in M$ with a representative $uv$ (resp. $vu$) for some $v \in A^\ast$, where for this final subset, we can take $L = uA^\ast$ (resp. $A^\ast u$). Therefore, decidability of the rational subset membership problem implies decidability of the word, submonoid membership, and divisibility problems for $M$. It is not difficult to show that, for \textit{any} monoid $M$ with context-free word problem, generated by a finite set $A$ and with associated surjective homomorphism $\pi \colon A^\ast \to M$, the pre-image $\pi^{-1}(R) \subseteq A^\ast$ of any rational subset $R \subseteq M$ is a context-free language; indeed, the proof is virtually identical to the proof of (i)$\implies$(ii) in Theorem~\ref{Thm:Equiv_of_three_notions}. Thus, from Theorem~\ref{Thm:Main_thm}, we deduce the following:

\begin{theorem}\label{Thm:RSMP_for_VF}
Let $M$ be a finitely presented special monoid with virtually free group of units. Then the rational subset membership problem for $M$ is decidable. 
\end{theorem}

This theorem generalises the fact that virtually free groups have decidable rational subset membership problem. Relatively much is known about the rational subset membership problem for groups, see e.g. \cite{Lohrey2008, Lohrey2015} and the survey by Lohrey \cite{Lohrey2015b}. This latter survey has 19 pages on the problem for groups, and yet a single paragraph (\S12) summarises all material known in 2015 regarding the problem for monoids. Thus Theorem~\ref{Thm:RSMP_for_VF} greatly expands the known results for monoids. For example, Render \& Kambites \cite{Render2009} prove that the bicyclic monoid $\pres{Mon}{b,c}{bc=1}$ has decidable rational subset membership problem. As the bicyclic monoid has trivial group of units, their result is a very special case of Theorem~\ref{Thm:RSMP_for_VF}.

\subsection{Other classes of languages}\label{Subsec:Other_classes}

In view of the restrictions on the class $\cc$ of languages to apply Theorem~\ref{Thm:Main_thm}, one might ask to what extent these may be weakened. The purpose of this section is to demonstrate that it is not at all straightforward to weaken them. First, we remark that the class $\REG$ of regular languages is not a super-$\AFL$, as it does not have the monadic ancestor property. The assumption of the monadic ancestor property cannot be directly removed from the statement of Theorem~\ref{Thm:Main_thm}, as in this case the bicyclic monoid (which has trivial, and hence regular) group of units does not have regular word problem, being infinite. However, it is not hard to understand the special monoids with regular word problem.

\begin{proposition}
Let $M$ be a finitely presented special monoid. Then $M$ has regular word problem if and only if $M$ is a finite group.
\end{proposition}
\begin{proof}
By An{\={\i}}s{\={\i}}mov's theorem, $M$ has regular word problem if and only if it is finite. On the other hand, from the results by Adian \cite{Adian1962} on identities in special monoids it follows that a special monoid is finite if and only if it is a group; see \cite[Theorem~6]{McNaughton1987} for a proof of this fact using rewriting techniques. 
\end{proof}

Another class of languages which is not a super-$\AFL$ is the class $\DCF$ of deterministic context-free languages. This is not, for example, closed under homomorphism or union. We can describe special monoids with deterministic context-free word problem quite well. First, free monoids have word problem in $\DCF$, as the language $\{ w \#w^\trev \mid w \in A^\ast \}$ is easily seen to be in $\DCF$. As an aside, this demonstrates the importance of the symbol $\#$, as the language of palindromes  $\{ w w^\trev \mid w \in A^\ast \}$ is not in $\DCF$ (see \cite[Exercise~12.6(b)]{Hopcroft1979}). 

We make an observation. It is not difficult to see that for any special monoid $M$ (generated by $A$ and with pieces $\Delta$) which is not isomorphic to a free product of a free monoid by a group, we have that $M$ contains a submonoid isomorphic to the bicyclic monoid. Indeed, $M$ is of the aforementioned form if and only if there is some piece $\delta \in \Delta$ with $|\delta|>1$. Write $\delta \equiv u_1 u_2$ with $u_1, u_2 \in A^+$, and let $u \in A^\ast$ be such that $u$ is an inverse of $\delta$. Let $v \equiv u_2 u$. Then it is not hard to see that while $u_1 v =_M 1$, we do not have $vu_1 =_M 1$. It follows by \cite[Lemma~1.31]{Clifford1961} that $\langle u_1, v \rangle \leq M$ is isomorphic to the bicyclic monoid. This is essentially the proof given in Lallement \cite[Theorem~1.2]{Lallement1974}, generalised from the one-relation case. Now, Brough, Cain \& Pfeiffer \cite{Brough2019} conjectured that the bicyclic monoid does not have deterministic context-free word problem; Kambites (unpublished) confirmed this conjecture, by an application of the pumping lemma for deterministic context-free languages \cite{Yu1989}. In particular, by \cite[Proposition~8(a)]{Hoffmann2012}, any monoid containing the bicyclic monoid does not have deterministic context-free word problem; by the above, thus any special monoid which is not a free product of a free monoid by a group does not have deterministic context-free word problem. It follows implicitly from Makanin \cite{Makanin1966b} or explicitly from Benois \cite{Benois1973} that a special monoid is isomorphic to a free product of a free monoid by a group if and only if it is right cancellative. We conclude:

\begin{proposition}\label{Prop:DCF=>RC+U(M)CF}
Let $M$ be a finitely presented special monoid with deterministic context-free word problem. Then $M$ is right cancellative, and the group of units $U(M)$ is a context-free group.
\end{proposition}

We emphasise that the case of having a special monoid be right cancellative is somewhat pathological. It is natural to conjecture that the converse of the above proposition holds, too.

\begin{conjecture}\label{Conj:DCFisRCVF}
Let $M$ be a finitely presented special monoid. Then $M$ has deterministic context-free word problem if and only if $M$ is right cancellative and the group of units $U(M)$ is a context-free group.
\end{conjecture}

The class of monoids with word problem in $\cc = \CF$ is closed under taking free products \cite[Theorem~6.2]{Brough2019}. Indeed, the author has shown that the same is true for any super-$\AFL$ $\cc$ \cite{NybergBrodda2020c}. On the other hand, it is an open problem whether the same is true for the class of monoids with word problem in $\DCF$. If this were true, then this would imply Conjecture~\ref{Conj:DCFisRCVF} holds, in view of Benois' result, the fact that free monoids have word problem in $\DCF$, and that the class of groups with word problem in $\CF$ coincides with the class of groups with word problem in $\DCF$ \cite{Muller1985}.

\subsection{Infix presentations}

Recall that in \S\ref{Subsec:Controlling_pres} it was proved that every special monoid admits such a presentation (with pieces $\Delta$) satisfying the small subpiece condition. A stronger condition than the small subpiece condition is the \textit{infix} condition, namely the condition that $\Delta$ be an infix code; i.e. no piece contains a piece as a proper subword. For example, let
\begin{equation*}
M = \pres{Mon}{a,b,c,d}{(ab)(cd)(ab)=1} \quad \text{resp.} \quad M' = \pres{Mon}{a,b}{(ab)(aabb)(ab)=1}.
\end{equation*}
These presentations have pieces $\Delta = \{ ab, cd\}$ resp. $\Delta' = \{ ab, aabb \}$. Thus $M$ is given by an infix presentation, but $M'$ is not. To the author it appears unlikely that $M'$ admits an infix presentation. However, at present, we have no direct means of proving this suspicion. This makes the following question natural.

\begin{question}\label{Quest:all_infix_free?}
Does every special monoid admit an infix presentation?
\end{question}

We conjecture that the answer to Question~\ref{Quest:all_infix_free?} is negative. The only general result we are able to show towards answering Question~\ref{Quest:all_infix_free?} is the following.

\begin{proposition}\label{Prop:UMtrivial=>Infix}
Let $M = \pres{Mon}{A}{w_1 = 1, w_2 = 1, \dots, w_k = 1}$. If the group of units of $M$ is trivial, then $M$ admits an infix presentation. Furthermore, an infix presentation for $M$ can be effectively computed from the given presentation for $M$. 
\end{proposition}
\begin{proof}
Let $\Delta = \{ \delta_1, \delta_2, \dots, \delta_n \}$ be the pieces of the presentation. As every $\delta_i$ is invertible, and $U(M) = 1$, we have $\delta_i =_M 1$ for every $1 \leq i \leq n$. Furthermore, if $w_\mu \equiv \delta_{i_1} \delta_{i_2} \cdots \delta_{i_\ell}$ with $1 \leq \mu \leq k$ and pieces $\delta_{i_j} \in \Delta$ and $1 \leq i_j \leq n$ for $1 \leq j \leq \ell$, then $w_\mu =_M 1$ follows from the set of relations $\{ \delta_i = 1 \mid 1 \leq i \leq n\}$. Thus the given presentation for $M$ is equivalent to the presentation
\begin{equation}\label{Eq:delta1=1,delta2=1,...}
M = \pres{Mon}{A}{\delta_1 = 1, \delta_2 = 1, \dots, \delta_n = 1}.
\end{equation}
It is obvious that $\Delta$ is the set of minimal invertible pieces for this presentation, too. Suppose the presentation \eqref{Eq:delta1=1,delta2=1,...} is not infix. Then there are pieces $\delta, \delta' \in \Delta$ such that $\delta \equiv h_1 \delta' h_2$ for some $h_1, h_2 \in A^+$. Add the relation $h_1 h_2 = 1$ to \eqref{Eq:delta1=1,delta2=1,...}, at which point the relation $\delta = 1$ (i.e. $h_1 \delta' h_2 = 1$) is redundant, following from $h_1h_2 = 1$ and $\delta' = 1$. The resulting presentation clearly has pieces $\Delta'$, where $\Delta' \subseteq \Delta$. We then repeat the process, given above to \eqref{Eq:delta1=1,delta2=1,...}, to the new presentation. As every piece in $\Delta$ contains at most finitely many occurrences of other pieces, this eventually terminates in an infix presentation. 

As $U(M) = 1$, the word problem is decidable for $M$ by Makanin's theorem, and so every step above is effective. 
\end{proof}

For example, consider the monoid
\[
M_9 = \pres{Mon}{a,b,c,d}{bcabcd=1, abcd=1, bc=1}.
\]
It is easy to check that $\Delta = \{ bc, abcd \}$, so this is not an infix presentation. However, $U(M_9) = 1$, so applying the treatment in the proof of Proposition~\ref{Prop:UMtrivial=>Infix} we find 
\[
M_9 \cong \pres{Mon}{a,b,c,d}{abcd=1, bc=1} \cong \pres{Mon}{a,b,c,d}{ad=1, bc=1}
\]
and the set of pieces of the final presentation is $\Delta' = \{ ad, bc \}$. Indeed, $M_9$ is simply the free product of two copies of the bicyclic monoid. In particular this final presentation is infix. This method is rather useful, and can likely be extended somewhat (say, to when $U(M)$ is finite?). Note, however, that Proposition~\ref{Prop:UMtrivial=>Infix} does nothing for $M'$ above, as $U(M') \cong \mathbb{Z}$. 

\section*{Acknowledgements}

\noindent Most of the research in this article was carried out while the author was a Ph.D. student at the University of East Anglia, and the main results appear in the author's Ph.D. thesis. The author wishes to thank his supervisor Robert D. Gray for valuable feedback (on several iterations of the present article), and Tara Brough for discussions pertaining to Conjecture~\ref{Conj:DCFisRCVF}. Finally, the author wishes to thank the Dame Kathleen Ollerenshaw Trust for funding his current research.

\bibliography{wp_special_monoids} 

\begin{thebibliography}{10}

\bibitem{Adian1960}
S.~I. Adian.
\newblock On the embeddability of semigroups in groups.
\newblock {\em Soviet Math. Dokl.}, 1:819--821, 1960.

\bibitem{Adian1962}
S.~I. Adian.
\newblock Identities in special semigroups.
\newblock {\em Dokl. Akad. Nauk SSSR}, 143:499--502, 1962.

\bibitem{Adian1966}
S.~I. Adian.
\newblock {\em Defining relations and algorithmic problems for groups and
  semigroups}.
\newblock Proceedings of the Steklov Institute of Mathematics, No. 85 (1966).
  1966.

\bibitem{Aho1968}
Alfred~V. Aho.
\newblock Indexed grammars---an extension of context-free grammars.
\newblock {\em J. Assoc. Comput. Mach.}, 15:647--671, 1968.

\bibitem{Anisimov1971}
A.~V. An{\={\i}}s{\={\i}}mov.
\newblock The group languages.
\newblock {\em Kibernetika (Kiev)}, (4):18--24, 1971.

\bibitem{Anisimov1975}
A.~V. An{\={\i}}s{\={\i}}mov.
\newblock Languages over free groups.
\newblock In {\em Mathematical foundations of computer science ({F}ourth
  {S}ympos., {M}ari\'{a}nsk\'{e} {L}\'{a}zn\v{e}, 1975)}, pages 167--171.
  Lecture Notes in Comput. Sci., Vol. 32. 1975.

\bibitem{Araujo2017}
V\'{\i}tor Ara\'{u}jo and Pedro~V. Silva.
\newblock Geometric characterizations of virtually free groups.
\newblock {\em J. Algebra Appl.}, 16(9):1750180, 13, 2017.

\bibitem{Benois1973}
Mich\`ele Benois.
\newblock Simplifiabilit\'{e} et plongement dans un groupe des mono\"{\i}des
  quotients d'un mono\"{\i}de libre par une congruence de {T}hue unitaire.
\newblock {\em C. R. Acad. Sci. Paris S\'{e}r. A-B}, 276:A665--A668, 1973.

\bibitem{Berstel1985}
Jean Berstel and Dominique Perrin.
\newblock {\em Theory of codes}, volume 117 of {\em Pure and Applied
  Mathematics}.
\newblock Academic Press, Inc., Orlando, FL, 1985.

\bibitem{Book1982}
Ronald~V. Book.
\newblock Confluent and other types of {T}hue systems.
\newblock {\em J. Assoc. Comput. Mach.}, 29(1):171--182, 1982.

\bibitem{Book1982b}
Ronald~V. Book, Matthias Jantzen, and Celia Wrathall.
\newblock Monadic {T}hue systems.
\newblock {\em Theoret. Comput. Sci.}, 19(3):231--251, 1982.

\bibitem{Book1993}
Ronald~V. Book and Friedrich Otto.
\newblock {\em String-rewriting systems}.
\newblock Texts and Monographs in Computer Science. Springer-Verlag, New York,
  1993.

\bibitem{Brough2019}
Tara Brough, Alan~J. Cain, and Markus Pfeiffer.
\newblock Context-free word problem semigroups.
\newblock In {\em Developments in language theory}, volume 11647 of {\em
  Lecture Notes in Comput. Sci.}, pages 292--305. Springer, Cham, 2019.

\bibitem{Cain2012}
Alan~J. Cain and Victor Maltcev.
\newblock Context-free rewriting systems and word-hyperbolic structures with
  uniqueness.
\newblock {\em Internat. J. Algebra Comput.}, 22(7), 2012.

\bibitem{Cain2014}
Alan~J. Cain and Victor Maltcev.
\newblock Hopfian and co-{H}opfian subsemigroups and extensions.
\newblock {\em Demonstr. Math.}, 47(4):791--804, 2014.

\bibitem{Cain2016}
Alan~J. Cain and Markus Pfeiffer.
\newblock Decision problems for word-hyperbolic semigroups.
\newblock {\em J. Algebra}, 465:287--321, 2016.

\bibitem{Campbell1995}
C.~M. Campbell, E.~F. Robertson, N.~Ru\v{s}kuc, and R.~M. Thomas.
\newblock Semigroup and group presentations.
\newblock {\em Bull. London Math. Soc.}, 27(1):46--50, 1995.

\bibitem{Clifford1961}
A.~H. Clifford and G.~B. Preston.
\newblock {\em The algebraic theory of semigroups. {V}ol. {I}}.
\newblock Mathematical Surveys, No. 7. American Mathematical Society,
  Providence, R.I., 1961.

\bibitem{Cochet1971}
Yves Cochet.
\newblock {\em Sur l'alg\'ebricit\'e des classes de certaines congruences
  d\'efinies sur le mono\"ide libre}.
\newblock PhD thesis, Universit\'e Rennes I, 1971.

\bibitem{Dehn1911}
M.~Dehn.
\newblock \"{U}ber unendliche diskontinuierliche {G}ruppen.
\newblock {\em Math. Ann.}, 71(1):116--144, 1911.

\bibitem{Diekert2013}
Volker Diekert and Armin Wei\ss.
\newblock Context-free groups and their structure trees.
\newblock {\em Internat. J. Algebra Comput.}, 23(3):611--642, 2013.

\bibitem{Duncan2004}
Andrew Duncan and Robert~H. Gilman.
\newblock Word hyperbolic semigroups.
\newblock {\em Math. Proc. Cambridge Philos. Soc.}, 136(3):513--524, 2004.

\bibitem{Dunwoody1985}
M.~J. Dunwoody.
\newblock The accessibility of finitely presented groups.
\newblock {\em Invent. Math.}, 81(3):449--457, 1985.

\bibitem{Eilenberg1969}
Samuel Eilenberg and M.~P. Sch\"{u}tzenberger.
\newblock Rational sets in commutative monoids.
\newblock {\em J. Algebra}, 13:173--191, 1969.

\bibitem{Engelfriet1985}
Joost Engelfriet.
\newblock Hierarchies of hyper-{AFL}s.
\newblock {\em J. Comput. System Sci.}, 30(1):86--115, 1985.

\bibitem{Fischer1972}
J.~Fischer, A.~Karrass, and D.~Solitar.
\newblock On one-relator groups having elements of finite order.
\newblock {\em Proc. Amer. Math. Soc.}, 33:297--301, 1972.

\bibitem{Garreta2019}
Albert Garreta and Robert~D. Gray.
\newblock Equations and first-order theory of one-relator and word-hyperbolic
  monoids, 2019.

\bibitem{Gilman1996}
Robert~H. Gilman.
\newblock Formal languages and infinite groups.
\newblock In {\em Geometric and computational perspectives on infinite groups
  ({M}inneapolis, {MN} and {N}ew {B}runswick, {NJ}, 1994)}, volume~25 of {\em
  DIMACS Ser. Discrete Math. Theoret. Comput. Sci.}, pages 27--51. Amer. Math.
  Soc., Providence, RI, 1996.

\bibitem{Gilman2007}
Robert~H. Gilman, Susan Hermiller, Derek~F. Holt, and Sarah Rees.
\newblock A characterisation of virtually free groups.
\newblock {\em Arch. Math. (Basel)}, 89(4):289--295, 2007.

\bibitem{Glushkov1961}
V.~M. Glu{\v{s}}kov.
\newblock Abstract theory of automata.
\newblock {\em Uspehi Mat. Nauk}, 16(5 (101)):3--62, 1961.

\bibitem{Gray2021}
Robert~D. Gray, Pedro~V. Silva, and N\'{o}ra Szak\'{a}cs.
\newblock Algorithmic properties of inverse monoids with hyperbolic and
  tree-like {S}ch\"utzenberger graphs.
\newblock {\em Pre-print}, 2021.
\newblock Available at arXiv:1912.00950.

\bibitem{Greibach1970}
Sheila~A. Greibach.
\newblock Full {${\rm AFLs}$} and nested iterated substitution.
\newblock {\em Information and Control}, 16:7--35, 1970.

\bibitem{Harrison1978}
Michael~A. Harrison.
\newblock {\em Introduction to formal language theory}.
\newblock Addison-Wesley Publishing Co., Reading, Mass., 1978.

\bibitem{Hoffmann2012}
Michael Hoffmann, Derek~F. Holt, Matthew~D. Owens, and Richard~M. Thomas.
\newblock Semigroups with a context-free word problem.
\newblock In {\em Developments in language theory}, volume 7410 of {\em Lecture
  Notes in Comput. Sci.}, pages 97--108. Springer, Heidelberg, 2012.

\bibitem{Hopcroft1979}
John~E. Hopcroft and Jeffrey~D. Ullman.
\newblock {\em Introduction to automata theory, languages, and computation}.
\newblock Addison-Wesley Publishing Co., Reading, Mass., 1979.
\newblock Addison-Wesley Series in Computer Science.

\bibitem{Jantzen1981}
M.~Jantzen.
\newblock On a special monoid with a single defining relation.
\newblock {\em Theoret. Comput. Sci.}, 16(1):61--73, 1981.

\bibitem{Jantzen1988}
Matthias Jantzen.
\newblock {\em Confluent string rewriting}, volume~14 of {\em EATCS Monographs
  on Theoretical Computer Science}.
\newblock Springer-Verlag, Berlin, 1988.

\bibitem{Kashintsev1978}
E.~V. Kashintsev.
\newblock On the word problem for special semigroups.
\newblock {\em Izv. Akad. Nauk SSSR Ser. Mat.}, 42(6):1401--1416, 1440, 1978.

\bibitem{Kashintsev1993}
E.~V. Kashintsev.
\newblock On the satisfiability of the conditions {$C'(\frac 13)$} and {$C(4)$}
  for special homogeneous semigroups with defining words-degrees.
\newblock {\em Mat. Zametki}, 54(3):40--47, 158, 1993.

\bibitem{Kobayashi2000}
Yuji Kobayashi.
\newblock Finite homotopy bases of one-relator monoids.
\newblock {\em J. Algebra}, 229(2):547--569, 2000.

\bibitem{Kral1970}
Jaroslav Kr{\'{a}}l.
\newblock A modification of a substitution theorem and some necessary and
  sufficient conditions for sets to be context-free.
\newblock {\em Math. Systems Theory}, 4:129--139, 1970.

\bibitem{Lallement1974}
G\'{e}rard Lallement.
\newblock On monoids presented by a single relation.
\newblock {\em J. Algebra}, 32:370--388, 1974.

\bibitem{LeGall2014}
Fran\c{c}ois Le~Gall.
\newblock Powers of tensors and fast matrix multiplication.
\newblock In {\em I{SSAC} 2014---{P}roceedings of the 39th {I}nternational
  {S}ymposium on {S}ymbolic and {A}lgebraic {C}omputation}, pages 296--303.
  ACM, New York, 2014.

\bibitem{Lohrey2015b}
Markus Lohrey.
\newblock The rational subset membership problem for groups: a survey.
\newblock In {\em Groups {S}t {A}ndrews 2013}, volume 422 of {\em London Math.
  Soc. Lecture Note Ser.}, pages 368--389. Cambridge Univ. Press, Cambridge,
  2015.

\bibitem{Lohrey2008}
Markus Lohrey and Benjamin Steinberg.
\newblock The submonoid and rational subset membership problems for graph
  groups.
\newblock {\em J. Algebra}, 320(2):728--755, 2008.

\bibitem{Lohrey2015}
Markus Lohrey, Benjamin Steinberg, and Georg Zetzsche.
\newblock Rational subsets and submonoids of wreath products.
\newblock {\em Inform. and Comput.}, 243:191--204, 2015.

\bibitem{Lyndon1977}
Roger~C. Lyndon and Paul~E. Schupp.
\newblock {\em Combinatorial group theory}.
\newblock Ergebnisse der Mathematik und ihrer Grenzgebiete, Band 89.
  Springer-Verlag, Berlin-New York, 1977.

\bibitem{Magnus1932}
W.~Magnus.
\newblock Das {I}dentit\"{a}tsproblem f\"{u}r {G}ruppen mit einer definierenden
  {R}elation.
\newblock {\em Math. Ann.}, 106(1):295--307, 1932.

\bibitem{Magnus1966}
Wilhelm Magnus, Abraham Karrass, and Donald Solitar.
\newblock {\em Combinatorial group theory: {P}resentations of groups in terms
  of generators and relations}.
\newblock Interscience Publishers [John Wiley \& Sons, Inc.], New
  York-London-Sydney, 1966.

\bibitem{Makanin1966b}
G.~S. Makanin.
\newblock {\em On the {I}dentity {P}roblem for {F}initely {P}resented Groups
  and {S}emigroups}.
\newblock PhD thesis, Steklov Mathematical Institute, Moscow, 1966.

\bibitem{Makanin1966}
G.~S. Makanin.
\newblock On the identity problem in finitely defined semigroups.
\newblock {\em Dokl. Akad. Nauk SSSR}, 171:285--287, 1966.

\bibitem{McNaughton1987}
R.~McNaughton and P.~Narendran.
\newblock Special monoids and special {T}hue systems.
\newblock {\em J. Algebra}, 108(1):248--255, 1987.

\bibitem{Muller1983}
David~E. Muller and Paul~E. Schupp.
\newblock Groups, the theory of ends, and context-free languages.
\newblock {\em J. Comput. System Sci.}, 26(3):295--310, 1983.

\bibitem{Muller1985}
David~E. Muller and Paul~E. Schupp.
\newblock The theory of ends, pushdown automata, and second-order logic.
\newblock {\em Theoret. Comput. Sci.}, 37(1):51--75, 1985.

\bibitem{Myasnikov2003}
Alexei~G. Myasnikov and Vladimir Shpilrain.
\newblock Automorphic orbits in free groups.
\newblock {\em J. Algebra}, 269(1):18--27, 2003.

\bibitem{Neumann1967}
B.~H. Neumann.
\newblock Some remarks on semigroup presentations.
\newblock {\em Canadian J. Math.}, 19:1018--1026, 1967.

\bibitem{Nivat1966}
Maurice Nivat.
\newblock \'{E}l\'{e}ments de la th\'{e}orie g\'{e}n\'{e}rale des codes.
\newblock In {\em Automata {T}heory}, pages 278--294. Academic Press, New York,
  1966.

\bibitem{Novikov1955}
P.~S. Novikov.
\newblock {\em Ob algoritmi\v{c}esko\u{\i} nerazre\v{s}imosti problemy
  to\v{z}destva slov v teorii grupp}.
\newblock Trudy Mat. Inst. im. Steklov. no. 44. Izdat. Akad. Nauk SSSR, Moscow,
  1955.

\bibitem{Thesis}
C.-F. Nyberg-Brodda.
\newblock {\em {The Word Problem and Combinatorial Methods for Groups and
  Semigroups}}.
\newblock PhD thesis, University of East Anglia, UK, 2021.

\bibitem{NybergBrodda2021f}
Carl-Fredrik Nyberg-Brodda.
\newblock The word problem for free products of monoids and semigroups.
\newblock {\em Pre-print}, 2021.
\newblock In preparation, exp. 2021.

\bibitem{NybergBrodda2021b}
Carl-Fredrik Nyberg-Brodda.
\newblock The word problem for one-relation monoids: a survey.
\newblock {\em Semigroup Forum}, 103(2):297--355, 2021.

\bibitem{NybergBrodda2020c}
Carl-Fredrik Nyberg-Brodda.
\newblock On the word problem for compressible monoids.
\newblock {\em Pre-print}, December 2020.
\newblock Available online at arXiv:2012.01402.

\bibitem{NybergBrodda2021c}
Carl-Fredrik Nyberg-Brodda.
\newblock {A translation of G. S. Makanin's 1966 Ph.D. thesis "On the Identity
  Problem for Finitely Presented Groups and Semigroups"}.
\newblock February 2021.
\newblock Available online at arXiv:2102.00745.

\bibitem{NybergBrodda2021d}
Carl-Fredrik Nyberg-Brodda.
\newblock Non-finitely generated maximal subgroups of context-free monoids.
\newblock {\em Pre-print}, July 2021.
\newblock Available online at arXiv:2107.12861.

\bibitem{Otto1991}
F.~Otto and L.~Zhang.
\newblock Decision problems for finite special string-rewriting systems that
  are confluent on some congruence class.
\newblock {\em Acta Inform.}, 28(5):477--510, 1991.

\bibitem{Perrin1984}
Dominique Perrin and Paul Schupp.
\newblock Sur les mono\"{\i}des \`a un relateur qui sont des groupes.
\newblock {\em Theoret. Comput. Sci.}, 33(2-3):331--334, 1984.

\bibitem{Render2009}
Elaine Render and Mark Kambites.
\newblock Rational subsets of polycyclic monoids and valence automata.
\newblock {\em Inform. and Comput.}, 207(11):1329--1339, 2009.

\bibitem{Shneerson1972a}
L.~M. Shneerson.
\newblock Identities in semigroups with one defining relation.
\newblock {\em Logic, Algebra, and Computational Mathematics, Ivanovo Pedag.
  Institute}, 1(1--2):139--156, 1972.

\bibitem{Shneerson1972b}
L.~M. Shneerson.
\newblock Identities in semigroups with one defining relation. {II}.
\newblock {\em Logic, Algebra, and Computational Mathematics, Ivanovo Pedag.
  Institute}, 1(3--4):112--124, 1972.

\bibitem{Stallings1968}
John~R. Stallings.
\newblock On torsion-free groups with infinitely many ends.
\newblock {\em Ann. of Math. (2)}, 88:312--334, 1968.

\bibitem{Thue1914}
Axel Thue.
\newblock Problem \"uber {V}er\"anderungen von {Z}eichenreihen nach gegebenen
  {R}egeln.
\newblock {\em Christiana {V}idenskaps-{S}elskabs {S}krifter, {I.} {M}ath.
  naturv. {K}lasse}, 10, 1914.

\bibitem{Tseitin1958}
G.~S. Tseitin.
\newblock An associative calculus with an insoluble problem of equivalence.
\newblock {\em Trudy Mat. Inst. Steklov.}, 52:172--189, 1958.

\bibitem{Valiant1975}
Leslie~G. Valiant.
\newblock General context-free recognition in less than cubic time.
\newblock {\em J. Comput. System Sci.}, 10:308--315, 1975.

\bibitem{Whitehead1936}
J.~H.~C. Whitehead.
\newblock On equivalent sets of elements in a free group.
\newblock {\em Ann. of Math. (2)}, 37(4):782--800, 1936.

\bibitem{Williams2012}
Virginia~Vassilevska Williams.
\newblock Multiplying matrices faster than {C}oppersmith-{W}inograd [extended
  abstract].
\newblock In {\em S{TOC}'12---{P}roceedings of the 2012 {ACM} {S}ymposium on
  {T}heory of {C}omputing}, pages 887--898. ACM, New York, 2012.

\bibitem{Yu1989}
Sheng Yu.
\newblock A pumping lemma for deterministic context-free languages.
\newblock {\em Inform. Process. Lett.}, 31(1):47--51, 1989.

\bibitem{Zhang1992}
Louxin Zhang.
\newblock Applying rewriting methods to special monoids.
\newblock {\em Math. Proc. Cambridge Philos. Soc.}, 112(3):495--505, 1992.

\bibitem{Zhang1992d}
Louxin Zhang.
\newblock Congruential languages specified by special string-rewriting systems.
\newblock In {\em Words, languages and combinatorics ({K}yoto, 1990)}, pages
  551--563. World Sci. Publ., River Edge, NJ, 1992.

\bibitem{Zhang1992b}
Louxin Zhang.
\newblock A short proof of a theorem of {A}djan.
\newblock {\em Proc. Amer. Math. Soc.}, 116(1):1--3, 1992.

\end{thebibliography}
\bibliographystyle{plain}

\end{document}